\documentclass[preprint,3p,times]{elsarticle}

\usepackage{amsmath}
\usepackage{array}
\usepackage{graphicx,color}
\usepackage{xcolor}
\usepackage{multicol}
\usepackage{multirow}
\usepackage[utf8]{inputenc}
\usepackage{amsmath,amsthm,amssymb}
\usepackage{amssymb}
\usepackage{algorithm}
\usepackage{algpseudocode}
\usepackage{booktabs}
\usepackage{amsfonts}
\usepackage{float}
\usepackage{subcaption}
\usepackage{comment}
\usepackage{todonotes}
\usepackage{tabularx}
\usepackage{cleveref}

\def\x{{\bf x}}

\def\nodes{{\cal V}}

\def\VHp{V_{H,p}}
\def\VHpG{\VHp^\Gamma}

\def\bv{{\bf v}}

\def\bu{{\bf u}}
\def\bv{{\bf v}}
\def\bh{{\bf h}}

\def\dsp{\displaystyle}

\definecolor{myred}{RGB}{200,0,0}
\newcommand{\red}[1]{\textcolor{myred}{#1}}

\def\O{\Omega}

\def\to{\rightarrow}

\def\H10{H^1_{ \partial \O \setminus \partial \O_S}}

\def\dsp{\displaystyle}

\def\kb{}

\def\red{}

\def\pou{P}

\def\ndof{m}
\makeindex             

\journal{Computers and Mathematics with Applications}
\begin{document}

\begin{frontmatter}


\title{Nonlinear Multilevel Solution Strategies for Diffusive Wave Flood Models in Perforated Domains}


\author[uca]{Miranda Boutilier\corref{cor1}}
\ead{miranda.boutilier@univ-cotedazur.fr}
\author[uca,inria]{Konstantin Brenner}
\ead{konstantin.brenner@univ-cotedazur.fr}
\author[tue]{Victorita Dolean}
\ead{v.dolean.maini@tue.nl}

 \cortext[cor1]{Corresponding Author}

\affiliation[uca]{organization={Université Côte d'Azur,  Laboratoire J.A. Dieudonné, CNRS UMR 7351},
            city={Nice},
            postcode={06108}, 
            country={France}}

\affiliation[inria]{organization={Team COFFEE, INRIA Centre at Université Côte d'Azur},
             city={Nice},
            postcode={06108},
            country={France}}

 \affiliation[tue]{organization={Eindhoven University of Technology, Department of Mathematics and Computer Science},
            city={Eindhoven},
            postcode={5600 MB},
  country={The Netherlands}}

\begin{abstract}
This article investigates the numerical solution of the Diffusive Wave equation posed on domains containing a large number of polygonal perforations, motivated by urban flood modeling. Such geometries induce strong multiscale effects driven by geometric complexity, which significantly challenge the robustness of standard nonlinear and linear solvers.
The work builds on a multiscale coarse space previously introduced by the authors for linear Poisson problems on perforated domains. This low-dimensional space, constructed on a coarse polygonal partition and spanned by locally discrete harmonic (Trefftz-type) basis functions, is shown to remain effective for the linearized Diffusive Wave problems arising within Newton iterations. This enables the construction of robust two-level preconditioners for the resulting sequence of linear systems.

Beyond linearization, the main focus of this work is on the effective solution of the fully nonlinear problem. We assess and combine several Schwarz-based nonlinear preconditioning strategies, including a two-level RASPEN method and a two-step nonlinear method, using the same multiscale coarse space to ensure scalability. While the individual components are drawn from the existing literature, their combination provides a robust and practical solution strategy for a challenging nonlinear problem posed on highly perforated domains. A systematic comparison of the methods and a discussion of algorithmic complexity are presented.
The proposed approaches are validated through numerical experiments, including a realistic test case based on topographical data from the city of Nice.

\end{abstract}
\end{frontmatter}

\section{Introduction}

\floatstyle{ruled}
\newfloat{method}{htbp}{loa}
\floatname{method}{Method}

\newfloat{procedure}{htbp}{lop}
\floatname{procedure}{Procedure}
\crefname{procedure}{Procedure}{Procedures}
\Crefname{procedure}{Procedure}{Procedures}

In this article, we investigate Schwarz-based domain decomposition and multiscale strategies for the numerical solution of stationary and time-dependent nonlinear partial differential equations on geometrically complex domains. Our main model problem is the Diffusive Wave (DW) equation \cite{diffwave1,diffwave},
\begin{equation}\label{dw_pde}
\partial_t u-\mathrm{div}\!\left( c_f\, h(u,z_b(\mathbf{x}))^\alpha \, \|\nabla u\|^{\gamma-1} \nabla u \right)=f,
\end{equation}
which is commonly used to model overland flows and can be derived from the Shallow Water equations by neglecting inertia terms \cite{diffwave}.
Here, $z_b(\mathbf{x})$ denotes the ground elevation (bathymetry), $u$ is the water surface elevation, and the water depth is
$
h(u,z_b(\mathbf{x}))=\max\!\big(u-z_b(\mathbf{x}),0\big).
$
The parameter $c_f$ is a friction coefficient, $1\le \alpha\le 2$, $\gamma\le 1$, and $f$ typically represents rainfall input.
Equation \eqref{dw_pde} is \emph{doubly nonlinear} (through $h(u,z_b)^\alpha$ and the gradient dependence) and may be \emph{degenerate}: the effective diffusion can vanish as $h\to 0$ and becomes singular as $\|\nabla u\|\to 0$.

Our motivation is the simulation of urban flooding, with a particular focus on the area of Nice, France.
Reliable forecasts can support mitigation decisions such as the placement of protective infrastructure (e.g., dams, dikes, drainage networks). A key modeling difficulty is that small-scale structural features (buildings, walls, fences) can strongly impact the timing and spatial extent of flooding, and their representation in the numerical model is therefore crucial \cite{YuLane1,YuLane2,Abily}. The numerical experiments in this paper rely on an infra-metric description of the urban environment provided by M\'{e}tropole Nice C\^{o}te d'Azur (MNCA) \cite{buildings}. From a hydraulic viewpoint, these structures can be considered essentially impervious and are thus modeled as perforations (holes) in the computational domain.

To formalize this setting, let $D\subset\mathbb{R}^2$ be an open, simply connected polygonal domain, and let $(\O_{S,k})_k$ be a finite family of mutually disjoint open polygonal subdomains of $D$ (the perforations), such that $\O := D\setminus\overline{\O_S}$ is connected, where $\O_S=\bigcup_k \O_{S,k}$.
We consider the DW model on $\O$:
\begin{equation}\label{eq:diffwave}
\left\{
\begin{array}{rll}
\dsp \partial_t u-\text{div}\!\left(\kappa(u,\nabla u)\nabla u\right) &=& 0 \qquad \mbox{in} \qquad \O\times(0,T],\\
\dsp u &=& g \qquad \mbox{on} \qquad (\partial\O\setminus\partial\O_S)\times(0,T],\\
\dsp \kappa(u,\nabla u)\,\frac{\partial u}{\partial \mathbf{n}} &=& 0 \qquad \mbox{on} \qquad (\partial\O\cap\partial\O_S)\times(0,T],
\end{array}
\right.
\end{equation}
supplemented with an initial condition, where
$\kappa(u,\nabla u)=h(u,z_b(\mathbf{x}))^\alpha\|\nabla u\|^{\gamma-1}$ and $\mathbf{n}$ is the outward unit normal.

The perforated geometry introduces fine spatial scales that must be resolved locally, while the nonlinear and possibly degenerate diffusion requires robust global nonlinear solvers. Domain decomposition offers a natural ``divide-and-conquer'' framework, based on a coarse partition of $\O$ into subdomains together with a fine triangulation resolving the perforations. However, scalable performance typically requires effective coarse (global) information transfer, and for nonlinear problems this raises additional algorithmic questions.

More {precisely,} we will focus on Schwarz-based methods for nonlinear problems which are often categorized into two classes: Newton--Krylov--Schwarz (NKS) and Schwarz--Newton--Krylov (SNK) methods. In NKS methods \cite{nkdd}, one applies Newton (or inexact Newton) to the non-linear equation $F(\bu)=0$ and uses Schwarz methods as \emph{linear} preconditioners for the Jacobian systems solved by Krylov methods. In contrast, SNK methods first construct a \emph{nonlinearly preconditioned} system $\mathcal{F}(\bu)=0$ that shares the solutions of $F(\bu)=0$, and then apply (inexact) Newton to $\mathcal{F}$; the nonlinear preconditioner is built from local nonlinear subproblems and can therefore localize stiff nonlinearities and mitigate stagnation phenomena in Newton iterations.

A prototypical SNK approach is ASPIN \cite{aspin}, which defines $\mathcal{F}$ so that nonlinearities are more balanced across the domain via additive local nonlinear corrections, followed by an inexact Newton solve. MSPIN \cite{mspin} provides a multiplicative counterpart, typically less parallel but often more contractive, analogously to the relation between additive and multiplicative Schwarz (or Jacobi and Gauss--Seidel) methods. RASPEN \cite{raspen} extends ASPIN by using restricted prolongations and a partition-of-unity construction, and yields a nonlinear fixed-point map for which an exact Jacobian of the preconditioned system can be derived, enabling exact Newton on $\mathcal{F}_{\mathrm{RASPEN}}(\bu)=0$. Related substructuring variants include SRASPEN \cite{sraspen} (nonoverlapping), DNPEN \cite{dnpen} (Dirichlet--Neumann-based), and adaptive strategies that selectively deactivate nonlinear preconditioning when not needed \cite{adaptivenonlin}.

Finally, as in the linear case, coarse (two-level) mechanisms are often essential for robustness of these nonlinear Schwarz variants. For RASPEN, a Full Approximation Scheme (FAS) coarse level was proposed in \cite{raspen}, with a multiplicative counterpart for MSPIN. Linear coarse corrections for ASPIN were proposed in \cite{linearcoarse}, and numerical evidence for robust two-level variants of ASPIN/RASPEN, including additive and multiplicative Galerkin-based coarse corrections and comparisons across coarse spaces, is reported in \cite{heinlein}.

\paragraph{Robustness through multiscale coarse spaces}
In this work, the term \emph{multiscale} does not refer to highly oscillatory or high-contrast material coefficients.
Instead, it denotes scale separation induced by \emph{geometric complexity}.
The presence of a large number of perforations fragments the domain and introduces multiple interacting geometric length scales.

Achieving robustness therefore requires coarse spaces capable of capturing global connectivity patterns that are not aligned with the underlying subdomain partition. 

The multiscale coarse spaces employed in this paper are specifically designed for this purpose.
By relying on energy {minimizing} 
MsFEM/Trefftz-type construction, they encode information associated with perforations enabling efficient information transfer across the domain.
This type of coarse correction is essential not only for scalability with respect to the number of subdomains (as in the classical use of coarse corrections), but also for robustness with respect to the geometric multiscale nature of the problem.

\paragraph{Main contributions}
Our contributions lie at the intersection of challenging application-driven modeling and a systematic assessment of nonlinear domain decomposition strategies.
Rather than proposing a new nonlinear framework, we focus on understanding how existing linear and nonlinear Schwarz-based methods behave when confronted with doubly nonlinear PDEs on highly perforated, multiscale geometries, and on identifying effective combinations of coarse spaces and nonlinear solvers.
Specifically, our main contributions are as follows:

\begin{itemize}
\item \textbf{Multiscale coarse spaces for nonlinear solves on perforated domains.}
Building on the MsFEM/Trefftz coarse space introduced in \cite{boutilier} for linear Poisson problems on perforated domains, we demonstrate that the same spectral multiscale construction remains effective for \emph{nonlinear} PDEs posed on highly perforated geometries, where robustness with respect to geometric scales is critical.

\item \textbf{Scalable two-level preconditioning for Newton linearizations.}
We incorporate this coarse space into two-level Restricted Additive Schwarz (RAS) preconditioners for the linear systems arising in Newton--Krylov and related linearized nonlinear solvers, enabling robust global communication across $N_s$ subdomains in geometrically multiscale settings.

\item \textbf{Coarse-enabled nonlinear domain decomposition strategies.}
We study several established nonlinear domain decomposition solvers, including one- and two-level RASPEN \cite{raspen} and the two-step NKS--RAS strategy \cite{cai2011}, and show how they can be combined with multiscale coarse corrections.
The emphasis is on \emph{effective coupling, robustness, and systematic assessment}, rather than on the introduction of new nonlinear algorithms.

\item \textbf{Acceleration and algorithmic variants in difficult nonlinear regimes.}
We investigate practical acceleration mechanisms, including Anderson acceleration \cite{onelevelrasanderson} applied to coarse-corrected nonlinear iterations, and identify regimes in which such enhancements yield tangible improvements over baseline nonlinear Schwarz methods.

\item \textbf{Extensive comparison and validation on doubly nonlinear multiscale models.}
We provide a detailed performance and complexity analysis, together with extensive numerical experiments on the doubly nonlinear Diffusive Wave equation posed on realistic, highly perforated urban geometries.
The methods are compared in terms of robustness with respect to the number of subdomains, geometric multiscale effects, and overall computational cost.
\end{itemize}

\paragraph{Outline of the paper}
The remainder of the paper is organized as follows.
\Cref{sec:disc} describes the finite-element/finite-volume hybrid discretization used for the Diffusive Wave equation.
\Cref{sec:coarsespace} \kb{we recall the} multiscale coarse space introduced in \cite{boutilier}, which underpins the two-level methods studied in this work.
\Cref{sec:nonlin} details the nonlinear domain decomposition strategies considered and discusses their computational cost and complexity.
\Cref{sec:numresults} reports numerical experiments on both academic and large-scale urban geometries, including the Porous Medium equation and the Diffusive Wave model.
Finally, \Cref{sec:conc} summarizes the main findings and outlines directions for future research.

\section{Discretization of the Diffusive Wave equation}\label{sec:disc}

We now describe the numerical discretization of the Diffusive Wave model
\eqref{eq:diffwave}.
For completeness and reproducibility, we briefly summarize the main ingredients of the discretization, while referring the reader to the cited literature for further details.

Let $T_f>0$ denote the final simulation time, and let
$0=t_0 < t_1 < \cdots < t_{N_T}=T_f$ be a partition of the time interval, with
$\Delta t_n = t_{n+1}-t_n$.
Given an approximation $u^n$ at time $t_n$, we consider the following semi-implicit time discretization of \eqref{eq:diffwave}: find $u^{n+1} \approx u(\cdot,t_{n+1})$ such that
\begin{equation}\label{eq:diffwave_semi-disc}
\left\{
\begin{array}{rll}
\dsp \frac{u^{n+1}-u^n}{\Delta t_n}
-\mathrm{div}\!\left( c_f \kappa(u^{n+1}, \nabla u^n)\nabla u^{n+1}\right) &=& 0
\quad \text{in } \O, \\[6pt]
\dsp u^{n+1} &=& g
\quad \text{on } \partial \O \setminus \partial \O_S, \\[6pt]
\dsp c_f \kappa(u^{n+1}, \nabla u^n)\,
\frac{\partial u^{n+1}}{\partial \mathbf{n}} &=& 0
\quad \text{on } \partial \O \cap \partial \O_S,
\end{array}
\right.
\end{equation}
with a given initial condition $u^0$.
This scheme treats the nonlinear diffusion coefficient semi-implicitly, which leads at each time step to a nonlinear elliptic problem.

\paragraph{Finite element discretization}
Let $\mathcal{T}$ be a triangulation of $\O$ conforming with the subdomain partition
$(\O_j)_{j=1}^{N_s}$.
\kb{We introduce the continuous, piecewise linear finite element space
\[
V_h
= \{ v \in C^0(\overline{\O}) \mid v|_T \in \mathbb{P}_1 \;\; \forall T\in\mathcal{T} \}
\]
and consider its subspace $V_{h,0}$ defined by
$$
V_{h,0} = \{ v_h \in V_h \,  | \, v_h |_{\partial \O \setminus \partial \O_S} = 0\}.
$$
}
Let $g_h\in V_h$ be a suitable approximation of the boundary data $g$, the finite element formulation of \eqref{eq:diffwave_semi-disc} reads as follows:
find $u_h^{n+1}\in V_h$ satisfying
$u_h^{n+1}=g_h$ on $\partial\O\setminus\partial\O_S$ such that
$$
\frac{1}{\Delta t_n}\int_\O (u_h^{n+1}-u_h^n)\,v_h\,\mathrm{d}\mathbf{x}
+\int_\O c_f \kappa(u_h^{n+1},\nabla u_h^n)
\nabla u_h^{n+1}\cdot\nabla v_h\,\mathrm{d}\mathbf{x}
=0,
\quad \forall v_h\in V_{h,0}.
$$

\paragraph{Nodal notation}
Let $(\mathbf{x}_i)_{i=1}^{\ndof_\O}$ denote the set of mesh nodes, with index set
$\mathcal{N}=\{1,\dots,\ndof_\O\}$, and let $(\eta_i)_{i\in\mathcal{N}}$ be the associated nodal basis functions.
We denote by $u_i^{n+1}$ the nodal value of $u_h^{n+1}$ at $\mathbf{x}_i$.
For each node $i\in\mathcal{N}$, we introduce the standard connectivity sets
\[
\mathcal{T}_i=\{T\in\mathcal{T}\mid \mathbf{x}_i\in\overline{T}\},
\qquad
\mathcal{N}_i=\bigcup_{T\in\mathcal{T}_i}\{\,\ell\mid \mathbf{x}_\ell\in\overline{T}\,\},
\]
and we denote by
$\mathcal{N}_D=\{i\in\mathcal{N}\mid \mathbf{x}_i\in\partial\O\setminus\partial\O_S\}$
the set of Dirichlet nodes.

\paragraph{FV--FE stabilization}
To enhance robustness and stability, in particular in the presence of degeneracy of the diffusion coefficient, we adopt a Finite Volume--Finite Element (FV--FE) discretization.
This consists of mass lumping in the accumulation term and upwinding in the diffusion term, and is also referred to as a control-volume finite element method
\cite{fem,control1,control2}.

\paragraph{Diffusion term}
Using the identity $\sum_{\ell\in\mathcal{N}_T}\nabla\eta_\ell=0$ on each element $T$,
the diffusion contribution associated with node $i\in\mathcal{N}$ can be rewritten as
\[
\int_\O c_f \kappa(u_h^{n+1},\nabla u_h^n)
\nabla u_h^{n+1}\cdot\nabla\eta_i\,\mathrm{d}\mathbf{x}
= \sum_{\ell\in\mathcal{N}_i} c_f (u_\ell^{n+1}-u_i^{n+1})
\sum_{T\in\mathcal{T}_i\cap\mathcal{T}_\ell}
\int_T \kappa(u_h^{n+1},\nabla u_h^n)
\nabla\eta_\ell\cdot\nabla\eta_i\,\mathrm{d}\mathbf{x}.
\]

To handle the degeneracy of
$\kappa(u,\boldsymbol{\xi})=h(u,z_b)^\alpha\|\boldsymbol{\xi}\|^{\gamma-1}$ when $h(u,z_b)=0$,
we introduce an upwind evaluation of the water depth.
We define
\[
\tau_{i\ell,T}^n
= -c_f\bigl|\nabla u_h^n|_T\bigr|^{1-\gamma}
\int_T \nabla\eta_\ell\cdot\nabla\eta_i\,\mathrm{d}\mathbf{x},
\qquad
\tau_{i\ell}^n=\sum_{T\in\mathcal{T}_i\cap\mathcal{T}_\ell}\tau_{i\ell,T}^n,
\]
and set the upstream water depth
\[
h_{i\ell}^{n+1}=
\begin{cases}
h(u_i^{n+1},z_b(\mathbf{x}_i)),
& \text{if } \tau_{i\ell}^n(u_i^{n+1}-u_\ell^{n+1})\ge 0,\\
h(u_\ell^{n+1},z_b(\mathbf{x}_\ell)),
& \text{otherwise}.
\end{cases}
\]
The diffusion term is then approximated by
\[
\sum_{\ell\in\mathcal{N}_i}
\tau_{i\ell}^n\bigl(h_{i\ell}^{n+1}\bigr)^\alpha
(u_i^{n+1}-u_\ell^{n+1}),
\qquad i\in\mathcal{N}\setminus\mathcal{N}_D.
\]

\paragraph{Accumulation term}
Let $M_{i\ell}=\int_\O\eta_i\eta_\ell\,\mathrm{d}\mathbf{x}$ denote the entries of the consistent mass matrix, and define the lumped mass
$m_i=\sum_{\ell=1}^{\ndof_\O}M_{i\ell}$.
The accumulation term is approximated by
\[
\frac{1}{\Delta t_n}\int_\O (u_h^{n+1}-u_h^n)\eta_i\,\mathrm{d}\mathbf{x}
\approx \frac{m_i}{\Delta t_n}(u_i^{n+1}-u_i^n),
\qquad i\in\mathcal{N}\setminus\mathcal{N}_D.
\]

\paragraph{Fully discrete scheme}
The resulting FV--FE discretization reads
\[
\begin{cases}
\dsp
\frac{m_i}{\Delta t_n}(u_i^{n+1}-u_i^n)
+\sum_{\ell\in\mathcal{N}_i}
\tau_{i\ell}^n\bigl(h_{i\ell}^{n+1}\bigr)^\alpha
(u_i^{n+1}-u_\ell^{n+1})
=0,
& i\in\mathcal{N}\setminus\mathcal{N}_D,\\[8pt]
u_i^{n+1}=g_h(\mathbf{x}_i),
& i\in\mathcal{N}_D.
\end{cases}
\]

\paragraph{Algebraic formulation}
Let $\bu^{n+1}=(u_i^{n+1})_{i\in\mathcal{N}\setminus\mathcal{N}_D}$.
The nonlinear system at time step $n+1$ can be written compactly as
\begin{equation}\label{dweq}
F(\bu^{n+1})=0,
\end{equation}
where the $i$th component of $F$ is given by
\begin{equation}\label{eq:dweresidual}
F_i(\bu^{n+1})
=\frac{m_i}{\Delta t_n}(u_i^{n+1}-u_i^n)
+\sum_{\ell\in\mathcal{N}_i}
\tau_{i\ell}^n\bigl(h_{i\ell}^{n+1}\bigr)^\alpha
(u_i^{n+1}-u_\ell^{n+1}),
\qquad i\in\mathcal{N}\setminus\mathcal{N}_D,
\end{equation}
with Dirichlet values enforced strongly on $\mathcal{N}_D$.

\section{\kb{Multiscale coarse space on perforated domains}}
\label{sec:coarsespace}

We briefly describe the multiscale coarse space introduced in \cite{boutilier}, which is used in all two-level algorithms considered in this paper.
This coarse space is of Trefftz type and consists of piecewise harmonic basis functions defined with respect to a polygonal partitioning of the domain. Let $\left( \kb{D}_j \right)_{j=1,\ldots,N_s}$ be a nonoverlapping polygonal partitioning of a background domain $D$,
and let $\O_j = \kb{D}_j \cap \O$ denote the induced partitioning of the perforated domain $\O$.

We refer to $\left( \O_j \right)_{j=1,\ldots,N_s}$ as the coarse mesh over $\O$. The nonoverlapping skeleton is defined as
\(
\Gamma
=
\bigcup_{j=1,\ldots,N_s}
\partial \O_j \setminus \partial \O_S.
\)
Let $\left( e_k \right)_{k=1,\ldots,N_e}$ be a partitioning of $\Gamma$ into open planar segments (coarse edges).
The set of coarse grid nodes is \kb{defined by}
\[
\nodes
=
\bigcup_{k=1,\ldots,N_e} \partial e_k.
\]

Let $H^1_\Delta(\O)$ denote the space of piecewise harmonic functions
\[
H^1_\Delta(\O)
=
\left\{
u \in H^1(\O)
\;\middle|\;
(u|_{\O_j}, v)_{H^1(\O_j)} = 0
\quad
\forall v \in \kb{H^1_{\partial \O_j\setminus \partial \O_S}(\O_j)},
\; \forall j = 1,\ldots,N_s
\right\},
\]
\kb{where $H^1_{\partial \O_j\setminus \partial \O_S}(\O_j)$ denote a subspace of $H^1(\O_j)$ composed of functions with vanishing traces on $\O_j\setminus \partial \O_S$. We note that 
that elements of $H^1_\Delta(\O)$ weakly satisfy homogeneous Neumann boundary conditions on $\partial \O \cap \partial \O_S$.
} 
We define the coarse trace space
\[
V_{H,p}^\Gamma
=
\left\{
v \in C^0(\overline{\Gamma})
\;\middle|\;
v|_{e_k} \in \mathbb{P}_p(e_k),
\; k=1,\ldots,N_e
\right\},
\]
and the Trefftz space
\[
V_{H,p}
=
\left\{
v \in H^1_\Delta(\O)
\;\middle|\;
v|_\Gamma \in V_{H,p}^\Gamma
\right\}.
\]
\begin{figure}
\centering
\includegraphics[height=6.cm]{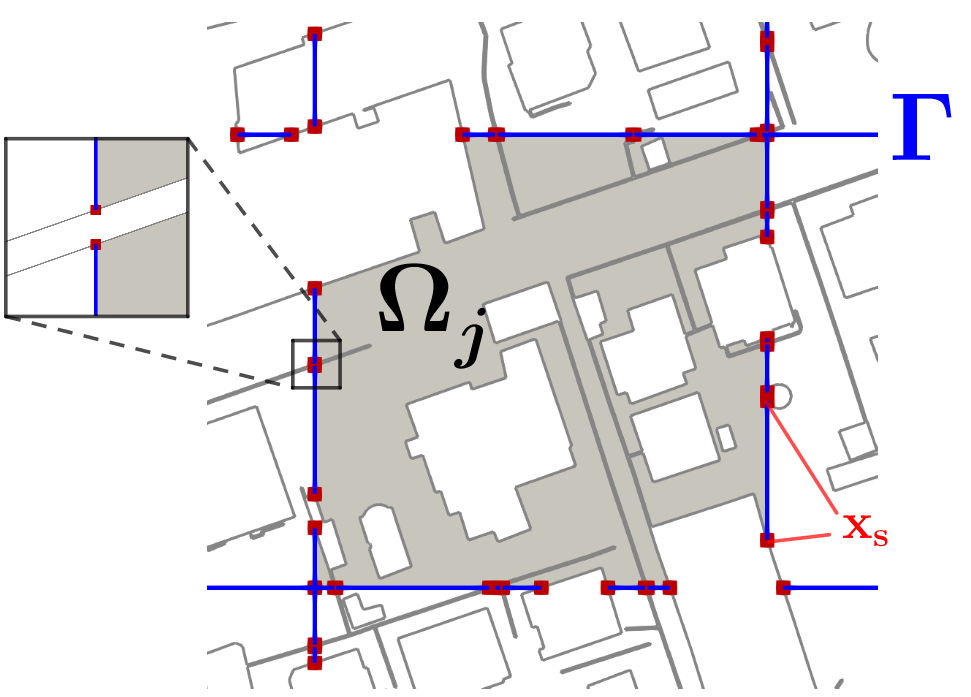}
\newline
\caption{Coarse cell $\O_j$, nonoverlapping skeleton $\Gamma$ (blue lines), and coarse grid nodes $\mathbf{x}_s \in \nodes$ (red dots). Coarse grid nodes are located at $\overline{\Gamma} \cap \partial \O_S.$}
\label{fig:coarsedofs}
\end{figure}
\kb{
Let $\left(g_s\right)_{s=1,\ldots,N_{\nodes}}$ denote the nodal Lagrange basis of $\VHpG$, and we will denote by $\left( \phi_s \right)_{s=1,\ldots,N_{\nodes}}$ the basis in $V_{H,p}$ obtained by the harmonic extension of $g_s$ over $\O$. That is, the restriction of $\phi_s$ to $\Omega_j$, denoted by $\phi_s^j$, weakly satisfies
$$
\left\{
\begin{array}{rll}
- \Delta \phi_s^j &=& 0 \quad \text{in } \O_j, \\
\dsp \frac{\partial \phi_s^j}{\partial \mathbf{n}} &=& 0 \quad \text{on } \partial \O_j \cap \partial \O_S, \\
\phi_s^j &=& g_s \quad \text{on } \partial \O_j \setminus \partial \O_S.
\end{array}
\right.
$$
}

Let $\overline{R}_j$ and $\overline{\pou}_j$ denote the restriction and partition-of-unity matrices associated with the nonoverlapping subdomains $\O_j$. \kb{Let $\boldsymbol{\phi}_s^j$ be a finite element approximation of $\phi_s^j$, t}he global coarse basis vectors are assembled as
\[
\boldsymbol{\phi}_s
=
\sum_{j \in \mathcal{N}_s}
\overline{R}_j^T \overline{\pou}_j \boldsymbol{\phi}_s^j,
\quad
s = 1,\ldots,N_{\nodes},
\]
where $\mathcal{N}_s$ is the set of subdomains containing the coarse node $\mathbf{x}_s$. The discrete Trefftz space is defined as
\(
\mathrm{span}\left( \left(\boldsymbol{\phi}_s\right)_{s=1,\ldots,N_{\nodes}}\right),
\)
and the coarse restriction matrix $R_H$ is formed by taking $\boldsymbol{\phi}_s^T$ as its rows.

\section{Nonlinear Multi-domain Solution Methods}
\label{sec:nonlin}
\kb{In this section, we present nonlinear solvers and preconditioning strategies for the discrete \kb{nonlinear system \eqref{dweq}.  Dropping the time step index, this system writes as $F(\bu) = 0$.} Let $\left( \O_j \right)_{j=1,\ldots,N_s}$ denote a nonoverlapping partitioning of $\O$, and let
$\left( \O_j' \right)_{j=1,\ldots,N_s}$ be an associated overlapping partitioning such that
$\O_j \subset \O_j'$.
In practice, each $\O_j'$ is obtained by extending $\O_j$ by a few layers of mesh elements.}

Let $R_j$ denote the Boolean restriction matrix from the global degrees of freedom to those associated with $\O_j'$, and let $\pou_j$ be partition-of-unity matrices satisfying
$$
\sum_{j=1}^{N_s} R_j^T \pou_j R_j = I.
$$

\kb{Below, we introduce several methods based on subdomain corrections using the overlapping decomposition $\left(\Omega_j'\right)_{j=1,\ldots,N_s}$. When a coarse second-level correction is included, it is constructed from the Trefftz coarse space described in Section~\ref{sec:coarsespace}, through the coarse restriction matrix $R_H$. The performance of these methods is reported in Section~\ref{sec:numresults}.}

\paragraph{Preconditioned Newton--Krylov method}
We first consider Newton's method applied to $F(\bu)=0$.
At outer iteration $k$, the Newton update $\boldsymbol{\delta}^k$ satisfies
\begin{equation}\label{eq:newtonlin}
\nabla F(\bu^k)\,\boldsymbol{\delta}^k = F(\bu^k),
\qquad
\bu^{k+1}=\bu^k-\boldsymbol{\delta}^k.
\end{equation}
The linear system \eqref{eq:newtonlin} is solved iteratively by GMRES preconditioned with the two-level RAS operator.
This gives the preconditioned system
$$
M_{\mathrm{RAS},2}^{-1}\nabla F(\bu^k)\,\boldsymbol{\delta}^k
=
M_{\mathrm{RAS},2}^{-1}F(\bu^k),
$$
\kb{where $M_{\mathrm{RAS},2}^{-1}$ is defined by}
\begin{equation}\label{eq:rasjn}
M_{\mathrm{RAS},2}^{-1}
=
\sum_{j=1}^{N_s} R_j^T \pou_j \left(R_j \nabla F(\bu^k) R_j^T\right)^{-1} R_j
+
R_H^T \left(R_H \nabla F(\bu^k) R_H^T\right)^{-1} R_H.
\end{equation}
\kb{In practice $M_{\mathrm{RAS},2}^{-1}$ is never explicitly assembled. Instead, the LU decomposition of $R_j \nabla F(\bu^k) R_j^T$ and  $R_H \nabla F(\bu^k) R_H^T$ is computed, which allows to define $M_{\mathrm{RAS},2}^{-1}$ as  a linear operator.}

\begin{algorithm}
\caption{Newton--Krylov with two-level RAS preconditioning to solve $F(\bu)=0$}
\label{alg:newtonmethod}
\begin{algorithmic}
\Require $\bu^0$
\For{$k=0,1,2,\ldots$ until convergence}
\State \kb{Assemble $F(\bu^k)$ and $\nabla F(\bu^k)$}
\State \kb{Setup the linear operator $M_{\mathrm{RAS},2}^{-1}$ from  \eqref{eq:rasjn}}
\State Approximately solve $\nabla F(\bu^k)\boldsymbol{\delta}^k = F(\bu^k)$ with GMRES preconditioned by $M_{\mathrm{RAS},2}^{-1}$
\State Set $\bu^{k+1} = \bu^k-\boldsymbol{\delta}^k$ %
\EndFor
\end{algorithmic}
\end{algorithm}

\paragraph{Nonlinear RAS (NRAS) and one-level RASPEN}
We now introduce the nonlinear analogue of the one-level RAS iteration.
For a given global vector $\bu$, define the local nonlinear \emph{subdomain update} $G_j(\bu)$
as the solution of
\begin{equation}\label{eq:subsolve}
R_j F\!\left(R_j^T G_j(\bu) + (I-R_j^TR_j)\bu\right)=0.
\end{equation}
This corresponds to solving the nonlinear problem on $\O_j'$ with Dirichlet data taken from $\bu$ on \kb{$\partial \O_j'\setminus \partial \O_S$}. 
In practice, \eqref{eq:subsolve} is solved, on each subdomain, by 
Newton's method; the associated local linear systems
are handled with a direct solver. The NRAS mapping is defined by assembling the local updates with the partition of unity:
\begin{equation}\label{eq:nras}
\mathrm{NRAS}(\bu) = \sum_{j=1}^{N_s} R_j^T \pou_j\, G_j(\bu).
\end{equation}
This yields the NRAS fixed-point iteration
$$
\bu^{k+1} = \mathrm{NRAS}(\bu^k).
$$
Applying Newton's method to the fixed-point equation
\begin{equation}\label{eq:raspen_fp}
\mathcal{F}(\bu) := \bu-\mathrm{NRAS}(\bu)=0
\end{equation}
yields the one-level RASPEN method \cite{raspen} (see also \cite{aspin} for related formulations).
We remark that the Jacobian of $\mathcal{F}$ is expected to be dense.  Specifically, we have
$$
\nabla \mathrm{NRAS}(\bu) = \sum_{j=1}^{N_s} R_j^T \pou_j\nabla G_j(\bu)
\qquad \mbox{with}\qquad
\nabla G_j(\bu)
=
-
\left( R_j \nabla F(\widetilde{\bu}_j) R_j^T \right)^{-1}
R_j \nabla F(\widetilde{\bu}_j) (I - R_j^T R_j),
$$
where
$$
\widetilde{\bu}_j = 
R_j^T G_j(\bu) + (I - R_j^T R_j)\bu.
$$
In practice $\nabla \mathcal{F}(\bu)$ is never explicitly assembled.  Instead,  the linear operator $\nabla \mathcal{F}(\bu)$ is setup based on the LU decomposition of $R_j J(\widetilde{\bu}_j) R_j^T$.
\begin{algorithm}
\caption{One-level RASPEN to solve $F(\bu)=0$}
\label{alg:1raspen}
\begin{algorithmic}
\Require $\bu^0$ 
\For{$k = 0,1,2,\ldots$ until convergence}
	\State \kb{For all $j$ compute $G_j(\bu^k)$ using Newton's method with direct linear solver}
    \State \kb{Set $\mathcal{F}(\bu^k) = \bu^k -  \mathrm{NRAS}(\bu^k)$ and setup the linear operator
    $\nabla \mathcal{F}(\bu^k) = I - \nabla \mathrm{NRAS}(\bu^k)$}
    \State Approximately solve 
    $\nabla \mathcal{F}(\bu^k)\boldsymbol{\delta}^k 
    = \mathcal{F}(\bu^k)$ with GMRES
    \State Set $\bu^{k+1} = \bu^k - \boldsymbol{\delta}^k$
\EndFor
\end{algorithmic}
\end{algorithm}

\paragraph{Two-level RASPEN}
One-level domain decomposition methods are not robust with respect to the number of subdomains $N_s$.
In particular, the number of Krylov iterations inside each outer iteration of Algorithm~\ref{alg:1raspen}
typically increases with $N_s$, motivating a two-level variant. We use the multiplicative two-level RASPEN strategy described in \cite[Algorithm~11]{thesisraspen},
where the coarse correction is computed after the local nonlinear solves.
Given $\bv\in\mathbb{R}^{\ndof_\O}$, the coarse nonlinear correction $c_H(\bv)$ is defined by the coarse problem
\begin{equation}\label{eq:coarsenonlin}
R_H F\!\left(\bv-R_H^T c_H(\bv)\right)=0.
\end{equation}
The resulting nonlinearly preconditioned residual is
\begin{equation}\label{eq:raspen2_res}
\mathcal{F}(\bu)
=
\bu-\mathrm{NRAS}(\bu) + R_H^T c_H(\mathrm{NRAS}(\bu)),
\end{equation}
and Newton's method is applied to $\mathcal{F}(\bu)=0$.
\begin{algorithm}[H]
\caption{Two-level RASPEN to solve $F(\bu)=0$}
\label{alg:2raspen}
\begin{algorithmic}
\Require $\bu^0$ 
\For{$k=0,1,2,\ldots$ until convergence}
\State \kb{For all $j$ compute $G_j(\bu^k)$ using Newton's method with direct linear solver and set $\widehat{\bu}^k = \mathrm{NRAS}(\bu^k)$}
\State Solve the coarse nonlinear problem \eqref{eq:coarsenonlin} for $c_H(\widehat{\bu}^k)$ using Newton's method with direct linear solver
\State Set $\mathcal{F}(\bu^k) = \bu^k-\widehat{\bu}^k + R_H^T c_H(\widehat{\bu}^k)$ and setup the linear operator
    $\nabla \mathcal{F}(\bu^k) = I - \nabla \mathrm{NRAS}(\bu^k) +  R_H^T \nabla c_H(\widehat{\bu}^k)$
    \State Approximately solve 
    $\nabla \mathcal{F}(\bu^k)\boldsymbol{\delta}^k 
    = \mathcal{F}(\bu^k)$ with GMRES
    \State Set $\bu^{k+1} = \bu^k - \boldsymbol{\delta}^k$
\EndFor
\end{algorithmic}
\end{algorithm}

\paragraph{Two-step NKS--RAS method}
We next consider the two-step method proposed in \cite{cai2011} (NKS--RAS).
At each outer iteration, this method performs one NRAS update combined with a global Newton correction based at the NRAS iterate.
It often converges rapidly in practice, but (except in special cases such as diagonal nonlinearities \cite{brennermonotone})
it cannot generally be interpreted as a nonlinear preconditioner in the sense of \kb{\eqref{eq:raspen_fp} or \eqref{eq:coarsenonlin}--\eqref{eq:raspen2_res}}. 
The linearized system is solved using the same two-level RAS-preconditioned GMRES strategy as in Algorithm~\ref{alg:newtonmethod}.

\begin{algorithm}
\caption{Two-step (NKS--RAS) method to solve $F(\bu)=0$}
\label{alg:2step}
\begin{algorithmic}
\Require $\bu^0$, nonlinear residual $F$, Jacobian $\nabla F$
\For{$k=0,1,2,\ldots$ until convergence}
\State Compute NRAS update $\widehat{\bu}^k = \mathrm{NRAS}(\bu^k)$
\State Approximately solve $\nabla F(\widehat{\bu}^k)\boldsymbol{\delta}^k = F(\widehat{\bu}^k)$ with GMRES preconditioned by $M_{\mathrm{RAS},2}^{-1}$
\State Set $\bu^{k+1} = \widehat{\bu}^k - \boldsymbol{\delta}^k$
\EndFor
\end{algorithmic}
\end{algorithm}

\paragraph{Two-level NRAS with Anderson acceleration}
The two-step method (Algorithm~\ref{alg:2step}) performs a global fine-scale Newton correction at each iteration.
We now replace the second step of the method by a coarse linear update based on the Trefftz coarse space.

For $\bv\in\mathbb{R}^{\ndof_\O}$, define the coarse linear correction operator
\begin{equation}\label{eq:coarselinear}
c_{H,l}(\bv)
=
\bv
-
R_H^T
\left(R_H \nabla F(\bv) R_H^T\right)^{-1}
R_H F(\bv).
\end{equation}
This yields the fixed-point iteration
\begin{equation}\label{eq:coarse2step}
\bu^{k+1} = c_{H,l}\!\left(\mathrm{NRAS}(\bu^k)\right).
\end{equation}
Since \eqref{eq:coarse2step} replaces a fine-scale Newton correction by a coarse update, its convergence is typically slower than that of Algorithm~\ref{alg:2step}. 
To improve convergence speed, we apply Anderson acceleration \cite{onelevelrasanderson} to the fixed-point map
$P(\bu) := c_{H,l}(\mathrm{NRAS}(\bu))$.
Specifically, at iteration $k$ we form residuals $\bv_i=P(\bu^i)-\bu^i$ and compute mixing coefficients
by a small constrained least-squares problem over the last $m$ iterates.

\begin{algorithm}
\caption{Anderson-accelerated coarse NRAS method to solve $F(\bu)=0$}
\label{alg:andersoncoarse2step}
\begin{algorithmic}
\Require $\bu^0$, nonlinear residual $F$, Jacobian $\nabla F$, history size $m\ge 1$
\State Define $P(\bu) := c_{H,l}(\mathrm{NRAS}(\bu))$ with $c_{H,l}$ given by \eqref{eq:coarselinear}
\State Set $\bu^1 = P(\bu^0)$
\For{$k=1,2,\ldots$ until convergence}
\State Set $m_k=\min\{m,k\}$ and form $\bv_i=P(\bu^i)-\bu^i$
\State Compute $\alpha^k\in\mathbb{R}^{m_k+1}$ solving
\[
\min_{\alpha^k}\left\|\sum_{i=0}^{m_k}\alpha_i^k\,\bv_{k-m_k+i}\right\|_2
\quad\text{s.t.}\quad
\sum_{i=0}^{m_k}\alpha_i^k=1
\]
\State Update $\bu^{k+1} = \sum_{i=0}^{m_k}\alpha_i^k\,P(\bu^{k-m_k+i})$
\EndFor
\end{algorithmic}
\end{algorithm}

\paragraph{Summary of the compared algorithms}
Table~\ref{tab:nonlin_methods_summary} summarizes the algorithms considered in this section,
highlighting the type of global solve required at each outer iteration and the use of coarse corrections.

\begin{table}[h!]
\centering
\caption{Summary of nonlinear multi-domain algorithms compared in this work ($N_s$ subdomains).}
\label{tab:nonlin_methods_summary}
\begin{tabular}{p{3.5cm} p{3.0cm} p{3.5cm} p{3.3cm}}
\toprule
\textbf{Algorithm} & \textbf{Outer iteration} & \textbf{Global linear solves} & \textbf{Coarse ingredient} \\
\midrule
Newton--Krylov + 2L RAS (Alg.~\ref{alg:newtonmethod})
& Newton on $F(\bu)=0$
& GMRES on $\nabla F(\bu^k)$
& 2L RAS with \kb{$R_H$}\\

1L RASPEN (Alg.~\ref{alg:1raspen})
& Newton on \kb{\eqref{eq:raspen_fp}}
& Linear solves with $\nabla\mathcal{F}(\bu^k)$ (operator form)
& none \\

2L RASPEN (Alg.~\ref{alg:2raspen})
& Newton on \eqref{eq:raspen2_res}
& Linear solves with $\nabla\mathcal{F}(\bu^k)$ (operator form)
& nonlinear coarse correction 
\eqref{eq:coarsenonlin} with $R_H$ \\

Two-step NKS--RAS (Alg.~\ref{alg:2step})
& NRAS + Newton correction
& GMRES on $\nabla F(\widehat{\bu}^k)$
& 2L RAS with \kb{$R_H$}\\

Anderson-accel.\ coarse NRAS (Alg.~\ref{alg:andersoncoarse2step})
& fixed point on \eqref{eq:coarse2step} + mixing
& coarse linear solves on $R_H \nabla F(\cdot) R_H^T$
& coarse linear correction \eqref{eq:coarselinear} with $R_H$ \\
\bottomrule
\end{tabular}
\end{table}

\subsection{Complexity of Fine-Scale Methods}\label{sec:complexity}

We summarize the computational cost of the methods introduced above.
All methods except Newton--Krylov involve a nonlinear restricted additive Schwarz (NRAS) update.

\paragraph{Scope and limitations of the complexity model}
We emphasize that the complexity estimates provided in this section are not intended
to predict wall-clock time or to compare methods in terms of absolute computational cost,
particularly in a parallel setting where communication, load imbalance, and implementation
details play a dominant role.
Rather, the purpose of this analysis is to identify the relative cost of the main algorithmic
building blocks (local nonlinear solves, global Krylov iterations, and coarse corrections)
and to clarify how different nonlinear strategies trade per-iteration complexity against
a reduction in the number of outer and inner iterations.

\paragraph{Notation summary}
For clarity, the main quantities used in the complexity estimates are summarized in
Table~\ref{tab:complexity_notation}.

\begin{table}[h]
\centering
\caption{Notation used in the complexity estimates.}
\label{tab:complexity_notation}
\renewcommand{\arraystretch}{1.2}
\begin{tabular}{ll}
\hline
Symbol & Meaning \\ \hline
$N_s$ & Number of subdomains \\
$\red{\ndof_{\O}}$ & Total number of degrees of freedom \\
$\red{\ndof_{\ell}}$ & Average number of local degrees of freedom per subdomain \\
$\red{\ndof_c}$ & Number of coarse degrees of freedom \\
$k_{\mathrm{out}}$ & Number of outer (nonlinear) iterations \\
$k_{\mathrm{loc}}$ & Average number of local Newton iterations per NRAS update \\
$k_{\mathrm{gm}}$ & Average number of GMRES iterations per outer iteration \\
$k_c$ & Average number of coarse linear solves per outer iteration \\
$C_{\mathrm{NRAS}}$ & Cost of one NRAS update \\
$C_{\mathrm{LU},\ell}$ & Cost of a local sparse LU factorization \\
$C_{\mathrm{LU},c}$ & Cost of a coarse sparse LU factorization \\
$C_{\mathrm{ass},\ell}$ & Cost of local residual/Jacobian assembly \\
$C_{\mathrm{ass}}$ & Cost of global residual/Jacobian assembly \\
$C_{\mathrm{mv}}$ & Cost of one sparse matrix--vector product \\
$N_T$ & Number of time steps (time-dependent problems) \\
\hline
\end{tabular}
\end{table}

\paragraph{Cost of one NRAS update}
An NRAS update requires the solution of independent nonlinear problems on overlapping subdomains.
Each local nonlinear problem is solved by Newton's method, with one sparse linear system solved at each local iteration. Let $N_s$ denote the number of subdomains, $k_{\mathrm{loc}}$ the average number of local Newton iterations,
$C_{\mathrm{LU},\ell}$ the cost of one local sparse LU factorization, and
$C_{\mathrm{ass},\ell}$ the cost of assembling the local residual and Jacobian.
The cost of one NRAS update is
$$
C_{\mathrm{NRAS}}
=
N_s \, k_{\mathrm{loc}}
\bigl(
C_{\mathrm{LU},\ell}
+
C_{\mathrm{ass},\ell}
\bigr).
$$
This expression assumes a sequential implementation; in parallel, the cost is dominated by the slowest subdomain.

\paragraph{Cost of GMRES solves}
All methods except Anderson acceleration require the solution of a global linear system using GMRES. Let $k_{\mathrm{gm}}$ denote the average number of GMRES iterations per outer iteration,
$C_{\mathrm{mv}}$ the cost of one sparse matrix--vector product,
$C_{\mathrm{LU},c}$ the cost of a coarse LU factorization, and
$C_{\mathrm{ass}}$ the cost of assembling the global residual and Jacobian. For Newton and two-step methods, the cost per outer iteration is
$$
C_{\mathrm{GMRES}}^{\text{(Newton/Two-step)}}
=
N_s \, C_{\mathrm{LU},\ell}
+
C_{\mathrm{LU},c}
+
C_{\mathrm{ass}}
+
k_{\mathrm{gm}} \, C_{\mathrm{mv}} .
$$

For RASPEN variants, the local factorizations and assemblies are reused from the NRAS step, leading to
$$
C_{\mathrm{GMRES}}^{\text{(RASPEN)}}
=
k_{\mathrm{gm}} \, C_{\mathrm{mv}} .
$$

\paragraph{Coarse nonlinear solves}
In the two-level RASPEN method, an additional nonlinear problem is solved on the coarse space.
Let $k_c$ denote the average number of coarse linear solves per outer iteration and
$C_{\mathrm{solve},c}$ the cost of one such solve.

\paragraph{Per-iteration cost summary}
The computational cost per outer iteration for each method is summarized as follows:
\begin{align*}
\textbf{Newton:} \quad
& C_{\mathrm{GMRES}}^{\text{(Newton)}} , \\[2mm]
\textbf{Two-step:} \quad
& C_{\mathrm{NRAS}} + C_{\mathrm{GMRES}}^{\text{(Two-step)}} , \\[2mm]
\textbf{One-level RASPEN:} \quad
& C_{\mathrm{NRAS}} + k_{\mathrm{gm}} \, C_{\mathrm{mv}} , \\[2mm]
\textbf{Two-level RASPEN:} \quad
& C_{\mathrm{NRAS}} + k_{\mathrm{gm}} \, C_{\mathrm{mv}} + k_c \, C_{\mathrm{solve},c} , \\[2mm]
\textbf{Anderson:} \quad
& C_{\mathrm{NRAS}} + C_{\mathrm{ass}} + O(\red{\ndof_{\O}}),
\end{align*}
where the $O(\red{\ndof_{\O}})$ term in Anderson acceleration corresponds to the least-squares problem.

\paragraph{Asymptotic costs}
Let $\red{\ndof_{\O}}$ denote the total number of degrees of freedom, $\red{\ndof_{\ell}} \approx \red{\ndof_{\O}} / N_s$ the number of local degrees of freedom per subdomain,
and $\red{\ndof_c}$ the number of coarse degrees of freedom.
Assuming that the LU factorization of a sparse $m \times m$ matrix scales as $O(m^{3/2})$ \cite{lu}, we obtain
\begin{align*}
C_{\mathrm{LU},\ell} &= O(\red{\ndof_{\ell}}^{3/2}), \\
C_{\mathrm{LU},c} = C_{\mathrm{solve},c} &= O(\red{\ndof_c}^{3/2}), \\
C_{\mathrm{ass},\ell} &= O(\red{\ndof_{\ell}}), \\
C_{\mathrm{ass}} = C_{\mathrm{mv}} &= O(\red{\ndof_{\O}}).
\end{align*}

\paragraph{Total cost}
Let $k_{\mathrm{out}}$ denote the number of outer iterations.
The total cost of a method is
\[
C_{\mathrm{total}} = k_{\mathrm{out}} \times (\text{per-iteration cost}).
\]
For time-dependent problems, all costs are multiplied by the number of time steps $N_T$, and
$k_{\mathrm{out}}$ denotes the average number of outer iterations per time step.

\section{Numerical Results}\label{sec:numresults}

We now present numerical experiments for the proposed algorithms (Algorithms 1–5) across several test cases. The experiments are arranged in order of increasing difficulty; they
include a stationary Porous Medium equation on an L-shaped domain, the similar equation on a small realistic urban domain, and the \kb{full} Diffusive Wave model on a large realistic urban domain. 
{\paragraph{On performance metrics}
We emphasize that the numerical experiments in this work are primarily intended
to compare the algorithmic efficiency and robustness of the proposed nonlinear solvers,
rather than to provide absolute performance measurements.
Since the current implementation is not optimized or parallel, raw CPU timings would be highly implementation-dependent
and potentially misleading.
For this reason, our main performance indicators are iteration-based and complexity-oriented metrics,
including the number of outer nonlinear iterations, cumulative GMRES iterations,
and the number of local nonlinear solves.
}
For all numerical experiments, the overlap is set to $\frac{1}{20}\mathcal{H}_j$, where $\mathcal{H}_j= \max (x_{max,j}-x_{min,j}, y_{max,j}-y_{min,j})$ and $x_{min,j}$, $y_{min,j}$, $x_{max,j}$, $y_{max,j}$ denote the minimal and maximal $x$ and $y$ coordinates that are contained in $\O_j$. We do not explore the effect of overlap further. {We note that in all of the considered test cases the triangulation is locally refined to account for the geometry constraints, and therefore the number of layers in the overlap is not fixed a priori.}

\subsection{Porous Medium Equation on L-shaped Domain}\label{sec:porous}

As a first numerical example we use an L-shaped domain with a reentering corner.  The domain is defined by $D = (-1,1)^2$,   $\O_S = (0,1)^2$ and $\O = D \setminus \overline{\O_S}$ such that the model domain has a singularity/corner in the upper right quadrant of the domain.
For this model domain, we consider the problem governed by the Porous Medium equation \kb{\cite{vazquez2007porous}}, given by
		\begin{equation}\label{eq:pme}
			\left\{
			\begin{array}{rll}
		 u(x,y) \kb{-} \mathrm{div}( \nabla(u(x,y)^4) ) &=& 0 \qquad \mbox{in} \qquad \O, \\
			\dsp  u(x,y)&=& 1 \qquad \mbox{on}  \qquad  y=1, \\
   			\dsp  u(x,y)&=& 0 \qquad \mbox{on} \qquad  x=1, \\
            \dsp  \frac{\partial u(x,y)}{\partial \mathbf{n}} &=& 0 \qquad \mbox{on} \qquad  \Gamma_N, \\
			\end{array}
			\right.
	\end{equation}
where $\Gamma_N=  \{(x, y) \in \partial \O \,|\, x \neq 1 \text{ and } y \neq 1\}$  denotes the Neumann boundary. 
Consider a piecewise affine finite element discretization, and let 
$\bu$ denote a vector of interior node values of size $\ndof_{\O}$. 
We denote the nodal ``hat" basis functions $(\eta^{\ell})_{\ell=1}^{\ndof_{\O}}$.
We define the lumped finite-element mass matrix $M$ and finite element stiffness matrix $A$ such that $$M_{ii}= \sum_{\ell=1}^{\ndof_{\O}} \int_\O \eta_i \eta_\ell \,\, {\rm d} \mathbf{x}, \qquad \text{and} \qquad A_{i\ell}= \int_\O \nabla \eta_i \cdot \nabla \eta_\ell \,\, {\rm d} \mathbf{x}.$$
\kb{The mass-lumped finite element discretization of \eqref{eq:pme} yields} 
\begin{equation}\label{eq:pmdedisc}
  F(\bu) := M \bu +A \max( \bu, 0 )^4 \kb{= 0}.
\end{equation}

\begin{figure}[t]
    \centering
    \includegraphics[width=0.4\linewidth]{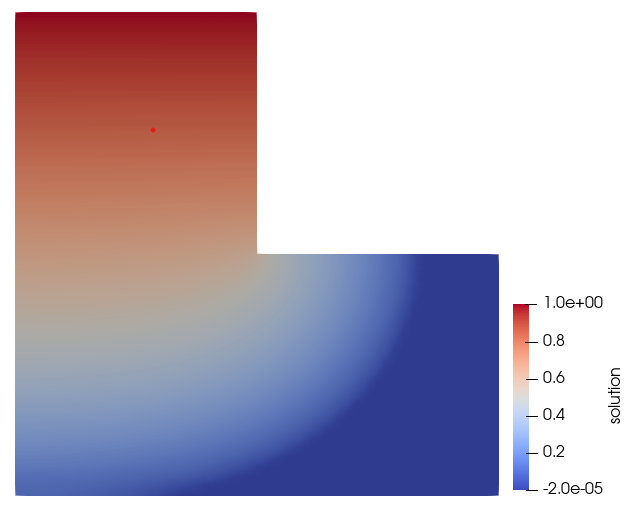}
    \caption{Finite element solution of \eqref{eq:pme} for test case \ref{sec:porous} on the chosen L-shaped domain.}
    \label{fig:pmesol}
\end{figure}

The finite element solution of this equation is shown in \Cref{fig:pmesol}. {In \eqref{eq:pmdedisc} the nonlinear term is approximated as
$
\left( \sum_{j} U_j \phi_j(x) \right)^4 \;\approx\; \sum_{j} U_j^4 \phi_j(x).
$
This nodal approximation is standard in finite volume / finite element and control volume finite element
discretizations of nonlinear diffusion 
equations \cite{control1,control2}.
Under Delaunay condition on the mesh, such formulations preserve positivity of the discrete solution and improve robustness in the presence of degenerate diffusion. 
}

\begin{figure}[t]
		\centering
		\begin{subfigure}{0.33\textwidth}
			\centering
			\includegraphics[width=\linewidth]{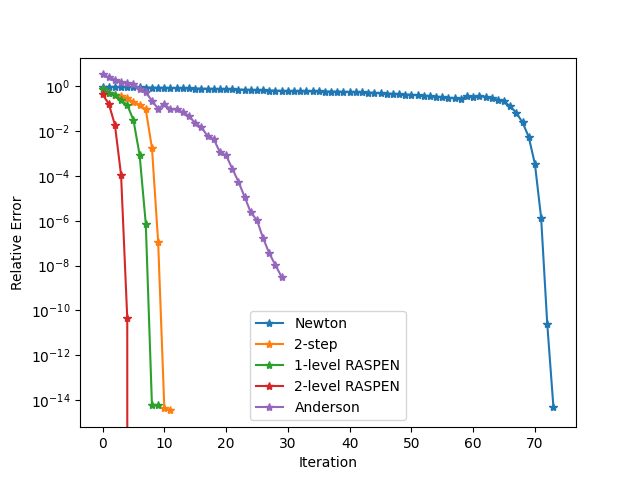}
		\end{subfigure}
		\begin{subfigure}{0.33\textwidth}
			\centering
			\includegraphics[width=\linewidth]{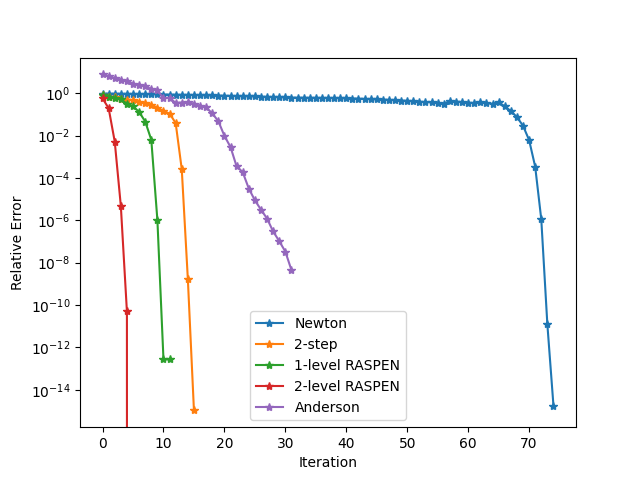}
		\end{subfigure}
  \begin{subfigure}{0.33\textwidth}
			\centering
			\includegraphics[width=\linewidth]{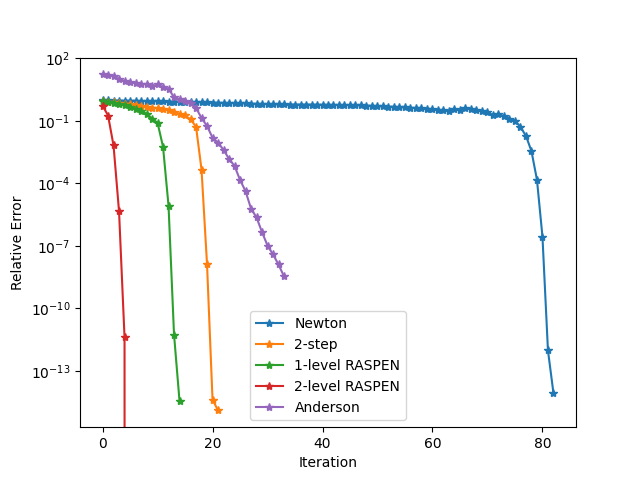}
		\end{subfigure}
		\caption{Test case \ref{sec:porous}, convergence curves. Top: $N=3$ (left), $N=5$ (right). Bottom: $N=9$.}
		\label{fig:convcurvespme}
	\end{figure}

\begin{table}[h]
\centering
\caption{Iteration counts for the test case \ref{sec:porous} for different numbers of subdomains $N_s$.
For two-level RASPEN, values in parentheses indicate $k_c$.}
\label{tab:ex2_compact}
\renewcommand{\arraystretch}{1.15}
\setlength{\tabcolsep}{3pt}
\begin{tabular}{lccc|ccc|ccc}
\hline
& \multicolumn{9}{c}{$N_s$} \\
Method
& \multicolumn{3}{c}{$3 \times 3$}
& \multicolumn{3}{c}{$5 \times 5$}
& \multicolumn{3}{c}{$9 \times 9$} \\
& $k_{\mathrm{out}}$
& $k_{\mathrm{gm}}$
& $k_{\mathrm{out}}k_{\mathrm{gm}}$
& $k_{\mathrm{out}}$
& $k_{\mathrm{gm}}$
& $k_{\mathrm{out}}k_{\mathrm{gm}}$
& $k_{\mathrm{out}}$
& $k_{\mathrm{gm}}$
& $k_{\mathrm{out}}k_{\mathrm{gm}}$ \\
\hline
Newton
& 79 & 33.7 & 2662
& 81 & 33.4 & 2705
& 88 & 38.1 & 3353 \\
Two-step
& 13 & 38.9 & 506
& 16 & 34.0 & 544
& 22 & 35.0 & 770 \\
1-level RASPEN
& 9 & 35.4 & 319
& 11 & 53.9 & 593
& 15 & 85.3 & 1280 \\
2-level RASPEN
& 6 & 20.5 (3.8) & 123 (23)
& 6 & 19.2 (4.5) & 115 (27)
& 6 & 19.2 (5.3) & 115 (32) \\
Anderson
& 25 & -- & --
& 32 & -- & --
& 31 & -- & -- \\
\hline
\end{tabular}
\end{table}

\paragraph{Methods and experimental setup}
We report numerical results for Newton's method \Cref{alg:newtonmethod},
the two-step method \Cref{alg:2step},
the one-level RASPEN method \Cref{alg:1raspen},
the two-level RASPEN method \Cref{alg:2raspen},
and Anderson acceleration \Cref{alg:andersoncoarse2step}.
For Newton and two-step methods, each outer iteration requires the solution of a linearized system
using GMRES preconditioned by the two-level linear RAS method \eqref{eq:rasjn}.

Results are presented for coarse partitionings
$N_s = 3 \times 3$, $5 \times 5$, and $9 \times 9$.
Choosing $N_s$ as the square of an odd integer ensures consistent mesh hierarchies
and comparable subdomain geometries across configurations.

\paragraph{Reported quantities}
Convergence histories are shown in \Cref{fig:convcurvespme}.
Iteration counts for all methods and all values of $N_s$ are summarized in
Table~\ref{tab:ex2_compact}.
The table reports the number of outer iterations $k_{\mathrm{out}}$,
the average number of GMRES iterations per outer iteration $k_{\mathrm{gm}}$ (when applicable),
and the cumulative number of GMRES iterations $k_{\mathrm{out}}\,k_{\mathrm{gm}}$.
For methods involving nonlinear subdomain solves, the table additionally reports
the average number of local linear solves per subdomain $k_{\mathrm{loc}}$.
For the two-level RASPEN method, the number of coarse linear solves $k_c$ is indicated in parentheses.

\paragraph{Outer iterations and convergence behavior}
From \Cref{fig:convcurvespme}, Newton's method exhibits a pronounced initial plateau
before entering the regime of quadratic convergence,
resulting in substantially larger values of $k_{\mathrm{out}}$ than all other methods.
In contrast, two-level RASPEN consistently achieves the smallest number of outer iterations
for all values of $N_s$.
The remaining methods rank, in increasing order of $k_{\mathrm{out}}$,
as one-level RASPEN, two-step, Anderson, and Newton.
While Anderson acceleration avoids the strong initial stagnation observed for Newton's method,
it does not attain the same asymptotic convergence rate as the nonlinear DD-based methods.

\paragraph{GMRES iterations and scalability}
Table~\ref{tab:ex2_compact} shows that two-level RASPEN is scalable
with respect to both $k_{\mathrm{out}}$ and $k_{\mathrm{gm}}$,
as these quantities remain approximately independent of $N_s$.
In contrast, the one-level RASPEN method exhibits a clear growth in $k_{\mathrm{gm}}$
as the number of subdomains increases, leading to a corresponding increase
in the cumulative GMRES count $k_{\mathrm{out}}\,k_{\mathrm{gm}}$.

Newton and two-step methods, which employ the same two-level linear RAS preconditioner,
yield comparable values of $k_{\mathrm{gm}}$ for all $N_s$.
However, because the two-step method converges in fewer outer iterations,
its cumulative GMRES cost $k_{\mathrm{out}}\,k_{\mathrm{gm}}$
is significantly lower than that of Newton's method.
Among all methods, two-level RASPEN achieves the smallest values of both
$k_{\mathrm{out}}$ and $k_{\mathrm{gm}}$,
and therefore the lowest total number of GMRES iterations.
Although two-level RASPEN additionally involves $k_{\mathrm{out}}\,k_c$ coarse linear solves,
these are inexpensive due to the small size of the coarse problem
and do not dominate the overall cost.

\paragraph{Local linear solves}
The average number of local linear solves per subdomain $k_{\mathrm{loc}}$
shows a mild dependence on $N_s$ for all methods,
which is expected since the local problem size
$\red{\ndof_{\ell}} \approx \red{\ndof_{\O}} / N_s$ increases as the number of subdomains decreases.
For a fixed $N_s$, the values of $k_{\mathrm{loc}}$ are comparable across all methods
that involve an NRAS update.
Consequently, differences in the total number of local linear solves per subdomain
are primarily driven by differences in the number of outer iterations $k_{\mathrm{out}}$.

\subsection{Porous Medium Equation on Large Urban Domain}\label{sec:medium}
\kb{As a second example we present numerical results for the Porous Medium equation on a large urban domain of size 
$160 \times160$ meters, $D = (-80, 80)^2$, containing a total of 72 perforations representing buildings and small fences. The geometry of the urban structures was provided by Métropole Nice Côte d’Azur.
}

\kb{For this example and the Diffusive Wave example below, we ensure that the fine-scale triangulation remains unchanged as the number of subdomains $N_s$ varies by first generating a sufficiently fine background mesh that can be used for all $N_s$. Specifically, we construct a fine-scale triangulation conforming to the coarse partition with $N_s = 16 \times 16$ subdomains and subsequently use this background triangulation for the cases $N_s = 2 \times 2$, $4 \times 4$, and $8 \times 8$ subdomains.  We use \cite{triangle} to generate the triangulation resulting in 35220 triangles and 21317 mesh points. Again, the overlap is set to $\frac{1}{20}\mathcal{H}_j$.}

\begin{figure}
\centering
\includegraphics[height=6cm]{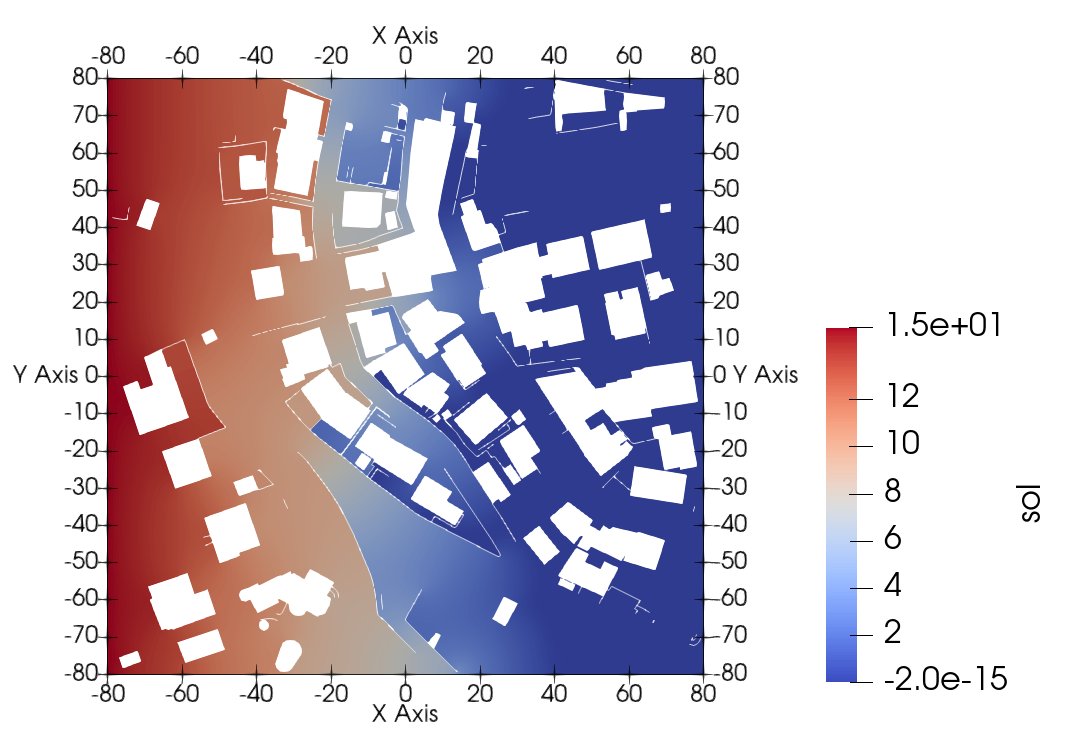}
\caption{Finite element solution for the test case \ref{sec:medium}.}
\label{fig:frameex3}
\end{figure}

We consider the Porous Medium equation given by
	\begin{equation}\label{eq:pmeex3}
			\left\{
			\begin{array}{rll}
		 u(x,y) \kb{-} c \, \mathrm{div}( \nabla(u(x,y)^3) ) &=& 0 \qquad \mbox{in} \qquad \O, \\
			\dsp  u(x,y)&=& c \qquad \mbox{on}  \qquad  x=-80, \\
   			\dsp  u(x,y)&=& 0 \qquad \mbox{on} \qquad  x=80, \\
            \dsp  \frac{\partial u(x,y)}{\partial \mathbf{n}} &=& 0 \qquad \mbox{on} \qquad  \Gamma_N, \\
			\end{array}
			\right.
	\end{equation}
\kb{where $\Gamma_N=  \{(x, y) \in \partial \O \,|\, x \neq -80 \text{ and } y \neq 80\}$  denotes the Neumann boundary,  which in particular includes the boundaries of the perforations. The parameter $c$ in \eqref{eq:pmeex3} is set to 15.}
With a similar finite element discretization to the test case \ref{sec:porous}, 
\kb{the discrete version of \eqref{eq:pmeex3} reads as
$$
F(\bu) :=  M \bu + c \,A \max(\bu,0)^m = 0,
$$
 The finite element solution for this example  shown in \Cref{fig:frameex3}.
}

\paragraph{Methods and \kb{experimental} setup}
We report numerical results for Newton's method \Cref{alg:newtonmethod},
the two-step method \Cref{alg:2step},
the two-level RASPEN method \Cref{alg:2raspen},
and Anderson acceleration \Cref{alg:andersoncoarse2step}.
Based on the results of the previous example, the one-level RASPEN method is omitted here,
as it is not robust for large numbers of subdomains.
For Newton and two-step methods, each outer iteration requires the solution of a linearized system
using GMRES preconditioned by the two-level linear RAS method \eqref{eq:rasjn}.

Results are shown for subdomain partitionings
$N_s = 2\times2$, $4\times4$, $8\times8$, and $16\times16$,
and for two initial guesses, $\mathbf u^0 = 1$ and $\mathbf u^0 = 0$.

\begin{figure}[t] 
		\centering
		\begin{subfigure}{0.45\textwidth}
		\centering
		\includegraphics[width=\linewidth]{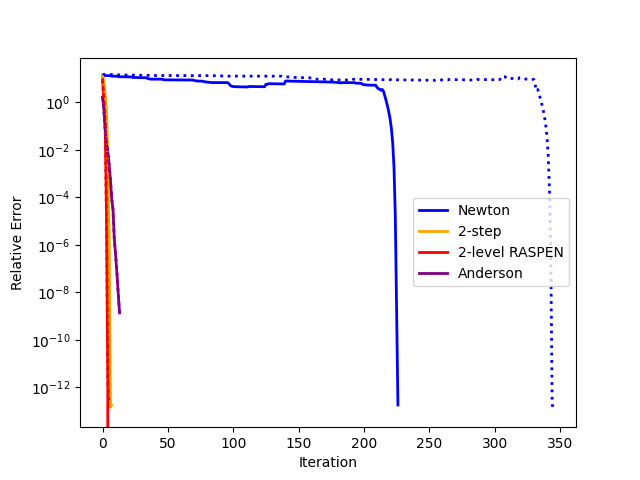}
	\end{subfigure}
	\begin{subfigure}{0.45\textwidth}
			\centering
			\includegraphics[width=\linewidth]{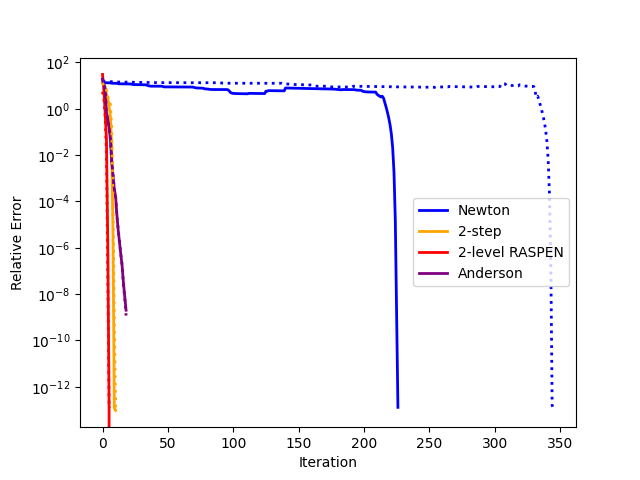}
		\end{subfigure}
  \begin{subfigure}{0.45\textwidth}
			\centering
			\includegraphics[width=\linewidth]{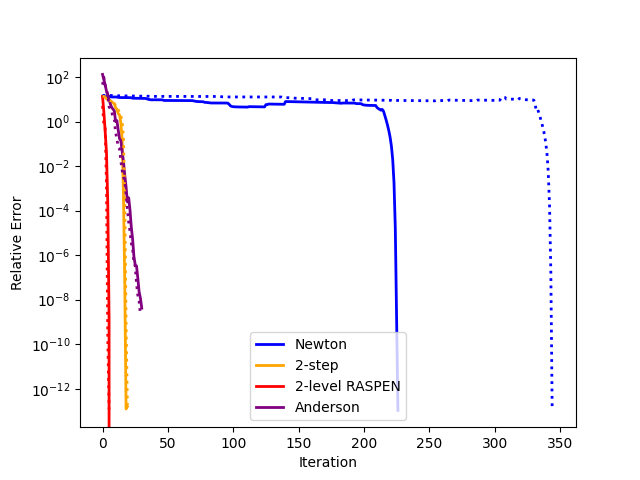}
		\end{subfigure}
    \begin{subfigure}{0.45\textwidth}
			\centering
			\includegraphics[width=\linewidth]{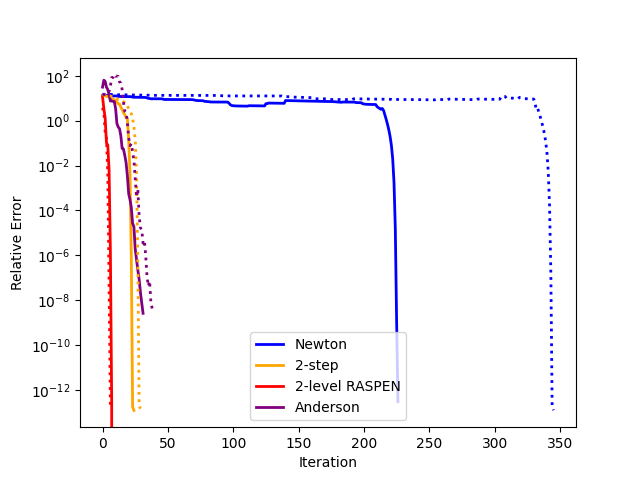}
		\end{subfigure}
		\caption{Test case \ref{sec:medium}, convergence curves. Solid lines correspond to an initial guess of $\bu^0=1$, dotted lines correspond to an initial guess of $\bu^0=0$. Top: $N_s=2 \times 2$ (left), $N_s=4 \times 4$ (right). Bottom: $N_s=8 \times 8$ (left), $N_s=16 \times 16$ (right).}
		\label{fig:convcurvesex3}
	\end{figure}

\paragraph{Reported quantities}
Convergence histories for both initial guesses are shown in \Cref{fig:convcurvesex3}.
Iteration counts are summarized in
Tables~\ref{tab:ex3_u1} and~\ref{tab:ex3_u0},
corresponding to the initial guesses $\mathbf u^0 = 1$ and $\mathbf u^0 = 0$, respectively.
These tables report the number of outer iterations $k_{\mathrm{out}}$,
the average number of GMRES iterations per outer iteration $k_{\mathrm{gm}}$ (when applicable),
and the cumulative number of GMRES iterations $k_{\mathrm{out}}k_{\mathrm{gm}}$.
For the two-level RASPEN method, the number of coarse linear solves $k_c$
is indicated in parentheses.

\paragraph{Outer iterations and robustness}
Newton's method exhibits a pronounced initial plateau before reaching the regime of quadratic convergence,
resulting in significantly larger values of $k_{\mathrm{out}}$ than all other methods.
This behavior is more pronounced than in the previous example and shows a strong dependence
on the initial guess (compare Tables~\ref{tab:ex3_u1} and~\ref{tab:ex3_u0}),
with increases of up to $50\%$ in $k_{\mathrm{out}}$ when switching from
$\mathbf u^0=1$ to $\mathbf u^0=0$.
In contrast, two-level RASPEN consistently achieves the smallest number of outer iterations
for all values of $N_s$ and for both initial guesses.
The two-step and Anderson methods exhibit a mild increase in $k_{\mathrm{out}}$ as $N_s$ increases,
but this growth is not proportional to $N_s$.

\paragraph{GMRES iterations and scalability}
Tables~\ref{tab:ex3_u1} and~\ref{tab:ex3_u0} show that the two-level RASPEN method
also yields the smallest values of $k_{\mathrm{gm}}$ and, consequently,
the lowest cumulative number of GMRES iterations $k_{\mathrm{out}}k_{\mathrm{gm}}$.
As in the previous example, Newton and two-step methods result in comparable values of
$k_{\mathrm{gm}}$ for all $N_s$.
However, the substantially larger number of outer iterations required by Newton's method
leads to between one and two orders of magnitude more total GMRES iterations.
The additional $k_{\mathrm{out}}k_c$ coarse linear solves required by two-level RASPEN
remain moderate and do not dominate the overall computational cost.

\paragraph{Summary}
Overall, the conclusions of this example are consistent with those of the test case \ref{sec:porous}.
Newton's method is particularly sensitive to the initial guess and suffers from a large number
of outer iterations, whereas two-level RASPEN remains robust with respect to both the number
of subdomains and the initial guess.
These results further confirm the effectiveness of the Trefftz coarse space
within a two-level nonlinear domain decomposition framework.

\begin{table}[h]
\centering
\caption{Iteration counts for the test case \ref{sec:medium} with initial guess $\mathbf u^0 = 0$
and different numbers of subdomains $N_s$.
For two-level RASPEN, values in parentheses indicate $k_c$.}
\label{tab:ex3_u0}
\renewcommand{\arraystretch}{1.15}
\setlength{\tabcolsep}{2pt}
\begin{tabular}{lccc|ccc|ccc|ccc}
\hline
& \multicolumn{12}{c}{$N_s$} \\
Method
& \multicolumn{3}{c}{$2\times2$}
& \multicolumn{3}{c}{$4\times4$}
& \multicolumn{3}{c}{$8\times8$}
& \multicolumn{3}{c}{$16\times16$} \\
& $k_{\mathrm{out}}$
& $k_{\mathrm{gm}}$
& $k_{\mathrm{out}}k_{\mathrm{gm}}$
& $k_{\mathrm{out}}$
& $k_{\mathrm{gm}}$
& $k_{\mathrm{out}}k_{\mathrm{gm}}$
& $k_{\mathrm{out}}$
& $k_{\mathrm{gm}}$
& $k_{\mathrm{out}}k_{\mathrm{gm}}$
& $k_{\mathrm{out}}$
& $k_{\mathrm{gm}}$
& $k_{\mathrm{out}}k_{\mathrm{gm}}$ \\
\hline
Newton
& 346 & 18.3 & 6332
& 346 & 25.4 & 8788
& 346 & 34.0 & 11764
& 346 & 46.0 & 15916 \\

Two-step
& 8 & 25.6 & 205
& 11 & 31.9 & 351
& 20 & 38.1 & 762
& 30 & 41.6 & 1248 \\

2-level RASPEN
& 6 & 13.2(4.3) & 79(26)
& 7 & 15.3(6.3) & 107(44)
& 6 & 19.5(8.5) & 117(51)
& 8 & 19.2(9.6) & 154(77) \\

Anderson
& 17 & -- & --
& 22 & -- & --
& 35 & -- & --
& 52 & -- & -- \\
\hline
\end{tabular}
\end{table}

\begin{table}[h]
\centering
\caption{Iteration counts for the test case \ref{sec:medium} with initial guess $\mathbf u^0 = 1$
and different numbers of subdomains $N_s$.
For two-level RASPEN, values in parentheses indicate $k_c$.}
\label{tab:ex3_u1}
\renewcommand{\arraystretch}{1.15}
\setlength{\tabcolsep}{2pt}
\begin{tabular}{lccc|ccc|ccc|ccc}
\hline
& \multicolumn{12}{c}{$N_s$} \\
Method
& \multicolumn{3}{c}{$2\times2$}
& \multicolumn{3}{c}{$4\times4$}
& \multicolumn{3}{c}{$8\times8$}
& \multicolumn{3}{c}{$16\times16$} \\
& $k_{\mathrm{out}}$
& $k_{\mathrm{gm}}$
& $k_{\mathrm{out}}k_{\mathrm{gm}}$
& $k_{\mathrm{out}}$
& $k_{\mathrm{gm}}$
& $k_{\mathrm{out}}k_{\mathrm{gm}}$
& $k_{\mathrm{out}}$
& $k_{\mathrm{gm}}$
& $k_{\mathrm{out}}k_{\mathrm{gm}}$
& $k_{\mathrm{out}}$
& $k_{\mathrm{gm}}$
& $k_{\mathrm{out}}k_{\mathrm{gm}}$ \\
\hline
Newton
& 227 & 21.0 & 4767
& 227 & 26.2 & 5947
& 227 & 35.2 & 7990
& 227 & 45.0 & 10215 \\

Two-step
& 8 & 25.6 & 205
& 11 & 30.1 & 331
& 20 & 38.4 & 768
& 25 & 41.1 & 1028 \\

2-level RASPEN
& 6 & 15.3(4.3) & 92(26)
& 7 & 20.4(5.9) & 143(41)
& 7 & 27.1(8.0) & 190(56)
& 9 & 22.3(7.4) & 201(67) \\

Anderson
& 17 & -- & --
& 22 & -- & --
& 35 & -- & --
& 37 & -- & -- \\
\hline
\end{tabular}
\end{table}

\subsection{Diffusive Wave Model on Large Realistic Urban Domain}\label{sec:LargeUrban_num}

In the final numerical experiment we consider the Diffusive Wave model \eqref{eq:diffwave} used to model a hypothetical flood in a densely urbanized area of the city of Nice, France.

\begin{figure} 
	\centering
		\begin{subfigure}{0.45\textwidth}
			\centering
			\includegraphics[width=.8\linewidth]{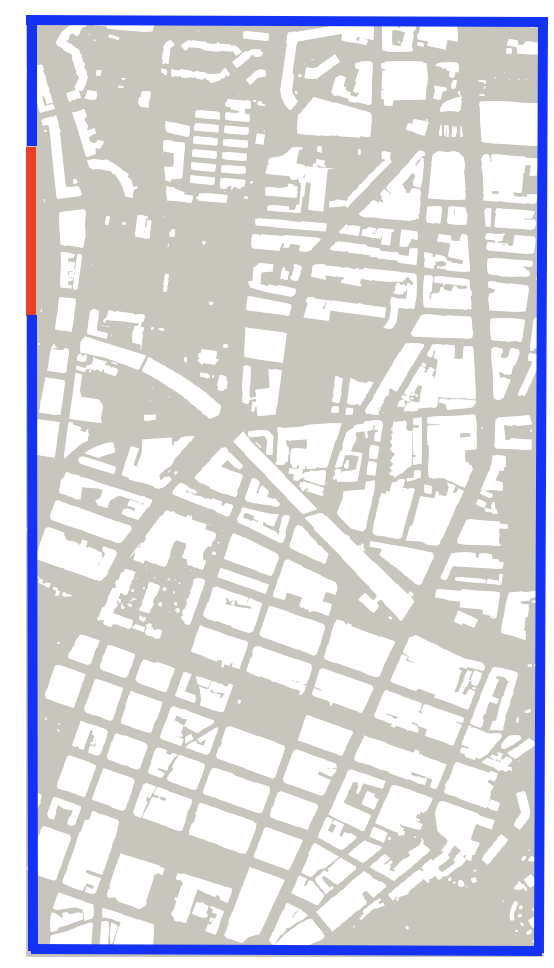}
		\end{subfigure}
  		\begin{subfigure}{0.45\textwidth}
			\centering
			\includegraphics[width=\linewidth]{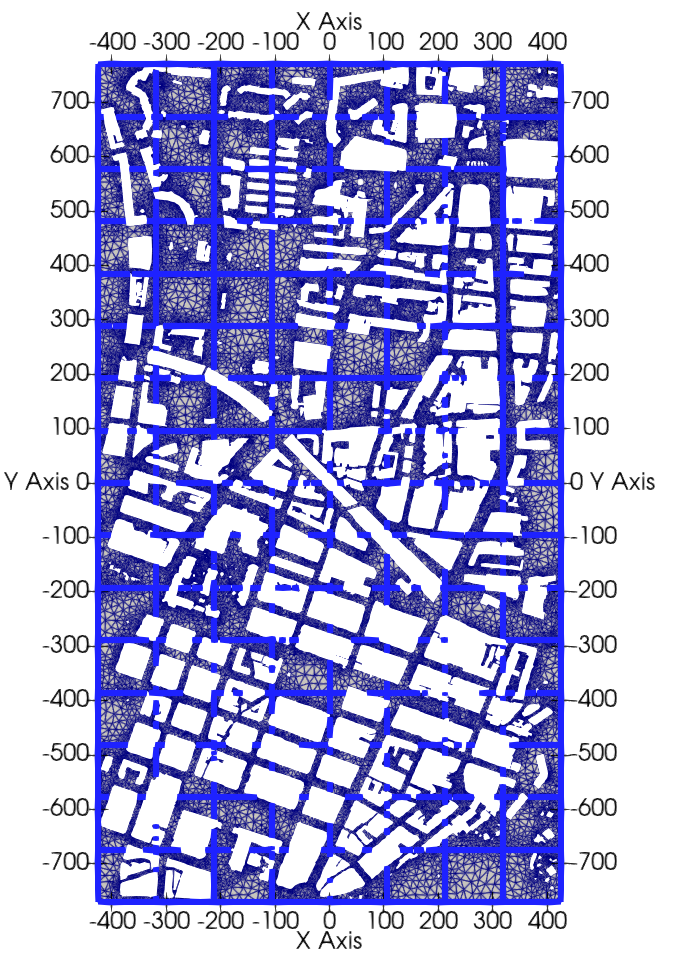}
		\end{subfigure}
        
		\caption{
        Left: Boundary condition (thick blue and red lines) and corresponding model domain for the test case \ref{sec:LargeUrban_num}. Red corresponds to $u=21.5$ meters, blue corresponds to $u=z_b$. This results in around a 2-meter elevation for the portion of the boundary representing the Paillon river. Right:
        Coarse (thick blue lines) and fine (thin blue lines) discretizations for $N_s=8 \times 16$ (right) subdomains.
}
		\label{fig:port}
	\end{figure}

\begin{figure} 
		\centering
		\begin{subfigure}{0.45\textwidth}
			\centering
			\includegraphics[width=\linewidth]{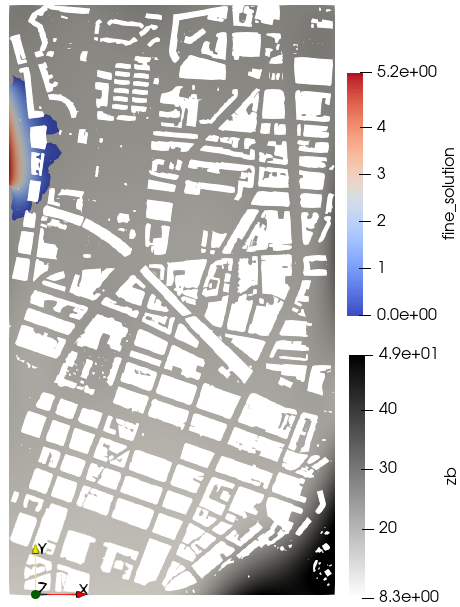}
		\end{subfigure}
		\begin{subfigure}{0.45\textwidth}
			\centering
			\includegraphics[width=\linewidth]{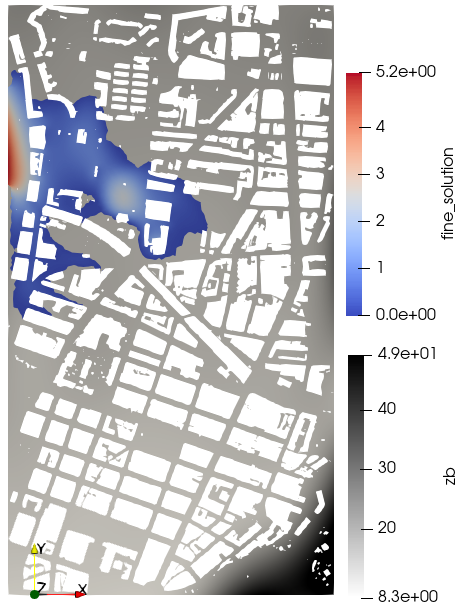}
		\end{subfigure}
  \newline 
  		\begin{subfigure}{0.45\textwidth}
			\centering
			\includegraphics[width=\linewidth]{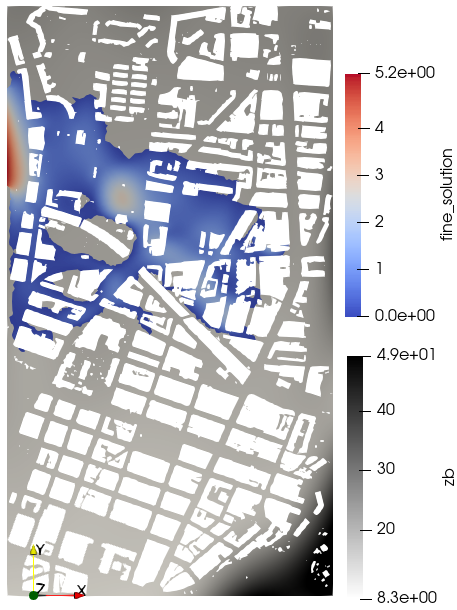}
		\end{subfigure}
		\begin{subfigure}{0.45\textwidth}
			\centering
			\includegraphics[width=\linewidth]{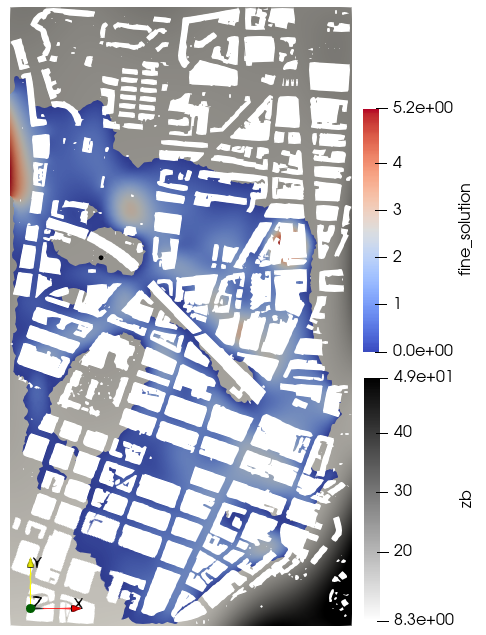}
		\end{subfigure}
		\caption{Numerical solution of for the test case \ref{sec:LargeUrban_num} at various time steps, where the numerical solution is presented by $\bh=\bu-z_b(\mathbf{x})$. Top: $t=10$ (left) seconds, $t=250$ (right) seconds. Bottom: $t=500$ (left) seconds, final $t=1500$ (right) seconds.  ``Flooded" region (where $\bh \geq 1$cm) are shown in color, with underlying bathymetry $z_b$ shown in the background in black and white. 
  }
		\label{fig:intialandfinaltime}
	\end{figure}

The model domain is depicted in \Cref{fig:intialandfinaltime} and has dimension $850 \times 1500$ meters.  In this test case scenario, the flood is produced by an overflow of the Paillon River in the north-west part of the domain. 
The buildings are removed from the computational domain, resulting in a total of 447 perforations. The geometry of the buildings was provided by Métropole Nice Côte d’Azur.
The bathymetry $z_b$ uses the $1m$ Digital Elevation Model available from \cite{elevationdata}.
The parameters of the Diffusive Wave model \eqref{eq:diffwave} are chosen as $\alpha=1.5$ and $\gamma=0.5$, based on Ch{\'e}zy's formula, with friction coefficient of $c_f=30$, which corresponding to a rough terrain \cite{white2011fluid}.
\kb{To model the overflow of the Paillon River, we prescribe a constant water elevation $u = 21.5$ meters above sea level on the corresponding portion of the boundary, which is approximately 2 meters above the riverbed, and we set $u=z_b$ elsewhere over the outer boundary. The Dirichlet boundary condition is shown in \Cref{fig:port}. The remainder of the boundary is subject to homogeneous Neumann boundary conditions.} The initial condition is taken to be equal to $z_b$, and the final simulation time is set to $T_f = 1500$ seconds. \kb{At this time, the solution has approximately reaches a steady state.}

Similarly to the test case \ref{sec:medium}, we keep the fine-scale triangulation fixed as the number of subdomains varies.
We first generate a background triangulation conforming to the finest partitioning $N_s=8\times16$,
and reuse this mesh for the coarser partitionings $N_s=1\times2$, $2\times4$, and $4\times8$.
\Cref{fig:port} illustrates the resulting consistent fine mesh for $N_s=2\times4$ and $N_s=8\times16$.

\paragraph{Local adaptive time stepping for NRAS}
At each time step, the nonlinear system \eqref{eq:dweresidual} must be solved.
For nonlinear DD methods, an NRAS update \eqref{eq:nras} requires the solution of local nonlinear subproblems.
In practice, 
a small subset of local subproblems
may fail to converge (typically due to stagnation).
To address this, we employ a local time-step reduction strategy:
failed subproblems are first solved with a smaller local time step and the resulting solution is used as an initial guess
for the original time step.  The procedure is summarized in Procedure \ref{alg:timestep}.

In this experiment, all methods except Newton's method use a fixed global time increment $\Delta t=10$\,s,
together with the local adaptive strategy of Procedure \ref{alg:timestep}.
We found that this avoids the need for global time-step reduction.
For Newton's method, we instead use a global adaptive strategy:
starting from $\Delta t_0=10$\,s, we reduce $\Delta t_n$ by a factor $\sqrt{2}$ upon failure,
and increase it by factor $\sqrt{2}$ after a successful solve, with $\Delta t_{n+1}=\min(\sqrt{2}\Delta t_n,10)$. The last time step is adapted to match the final time $T_f$.

\algnewcommand{\IfThenElse}[3]{
  \State \algorithmicif\ #1\ \algorithmicthen\ #2\ \algorithmicelse\ #3}
\algnewcommand{\IfThen}[2]{
  \State \algorithmicif\ #1\ \algorithmicthen\ #2}

\begin{procedure}
\caption{Continuation method for the parametrized system $F_{\Delta t}(\bu)=0$}
\label{alg:timestep}
\begin{algorithmic}
\Require Initial guess $\bu_0$

\State $\Delta t_{\text{worked}} = \Delta t$, $\bu = \bu_0$

\State Step 1: Find a working time step
\While{true}
    \State Attempt to compute $ \bu = F_{\Delta t_{\text{worked}}}^{-1}(0)$ using Newton's method with initial guess $\bu$
	\IfThenElse{Newton fails}{$\Delta t_{\text{worked}} = \Delta t_{\text{worked}}/2$}{\textbf{break}}
\EndWhile

\IfThen{$\Delta t_{\text{worked}} = \Delta t$}{\Return $\bu$}

\State Step 2: Progressively increase the time step
\While{true}

    \State Attempt to compute $\bu = F_{\Delta t}^{-1}(0)$ using Newton's method with initial guess $\bu$

	\IfThen{Newton succeeds}{\Return $\bu$}

    \State $\Delta t_{\text{try}} = (\Delta t_{\text{worked}}+\Delta t)/2$

    \While{true}
        \State Attempt to compute $\bu = F_{\Delta t_{\text{try}}}^{-1}(0)$ using Newton's method with initial guess $\bu$
 
		\IfThenElse{Newton fails}{$\Delta t_{\text{try}} = (\Delta t_{\text{worked}}+\Delta t_{\text{try}})/2$}{ $\Delta t_{\text{worked}} = \Delta t_{\text{try}}$ and \textbf{break}}
    \EndWhile

\EndWhile

\end{algorithmic}
\end{procedure}

\begin{figure}[t] 
		\centering
  	\begin{subfigure}{0.45\textwidth}
			\centering
			\includegraphics[width=\linewidth]{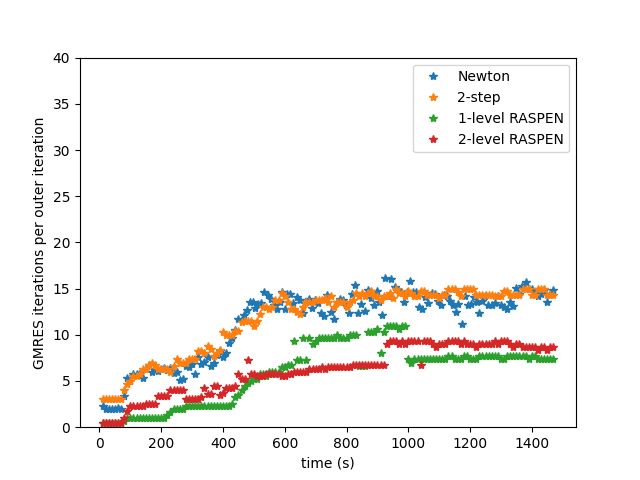}
		\end{subfigure}
		\begin{subfigure}{0.45\textwidth}
			\centering
			\includegraphics[width=\linewidth]{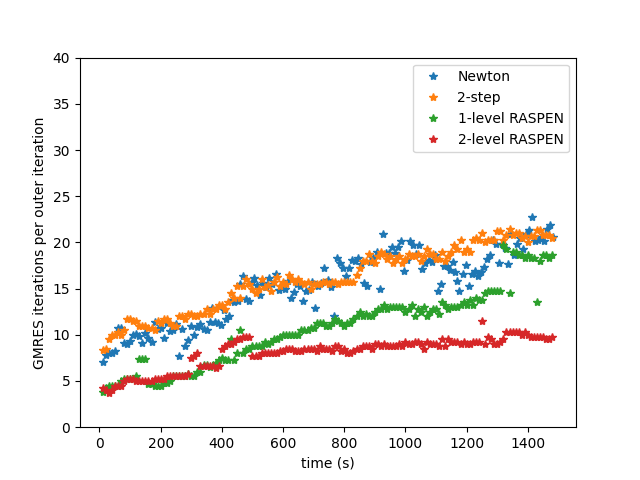}
		\end{subfigure}
		\begin{subfigure}{0.45\textwidth}
			\centering
			\includegraphics[width=\linewidth]{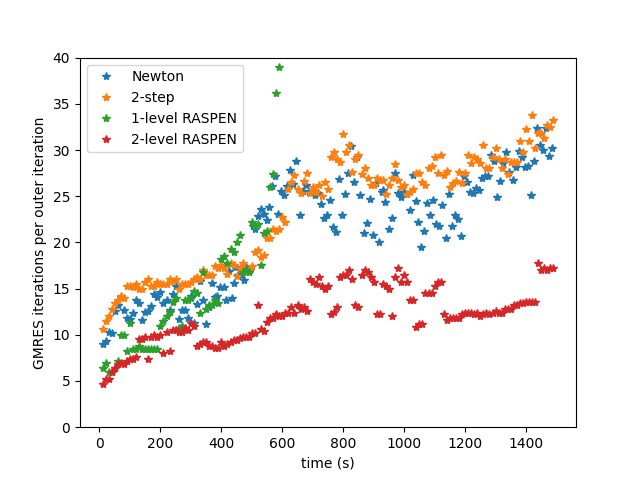}
		\end{subfigure}
  \begin{subfigure}{0.45\textwidth}
			\centering
			\includegraphics[width=\linewidth]{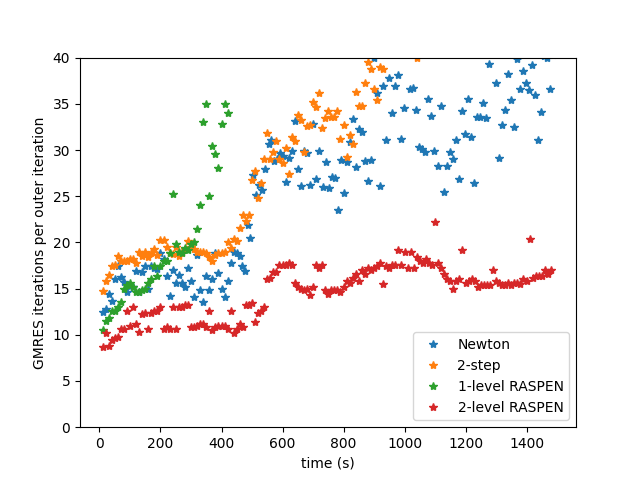}
		\end{subfigure}
		\caption{Test case \ref{sec:LargeUrban_num},  average number of GMRES iterations per outer iteration $k_{GM}$ at each time step. Top: $N_s=1 \times 2$ (left), $N_s=2\times 4$ (right). Bottom: $N_s=4 \times 8$ (left), $N_s=8 \times 16$ (right).}
 
		\label{fig:linsolvesperit}
	\end{figure}

\paragraph{Methods and reported quantities}
At each time step, we solve \eqref{eq:dweresidual} using
Newton's method (\Cref{alg:newtonmethod}),
the two-step method (\Cref{alg:2step}),
one- and two-level RASPEN (\Cref{alg:1raspen} and \Cref{alg:2raspen}),
and Anderson acceleration \Cref{alg:andersoncoarse2step}.
For Newton's and two-step methods, each outer iteration requires a GMRES solve preconditioned by
the two-level linear RAS preconditioner \eqref{eq:rasjn}.
The per-time-step GMRES iteration counts $k_{\mathrm{gm}}$ are shown in \Cref{fig:linsolvesperit}.
A summary of cumulative iteration counts over the full simulation is reported in
Table~\ref{tab:ex3_time_compact}.

\paragraph{Outer iterations}
From Table~\ref{tab:ex3_time_compact}, the two-step method and both RASPEN variants yield the smallest cumulative outer-iteration counts,
with similar behavior across all values of $N_s$.
Newton requires noticeably more outer iterations.
Anderson exhibits the strongest dependence on $N_s$, with cumulative outer iterations increasing as $N_s$ grows.

\paragraph{GMRES iterations and robustness with respect to $N_s$}
From \Cref{fig:linsolvesperit}, one-level RASPEN loses robustness as $N_s$ increases:
the per-iteration GMRES counts become substantially larger than for the other methods,
and for $N_s=8\times16$ they exceed the plotting range.
In contrast, two-level RASPEN yields the smallest $k_{\mathrm{gm}}$ as $N_s$ grows.
Newton and two-step methods, both preconditioned by the same two-level linear RAS operator,
exhibit comparable values of $k_{\mathrm{gm}}$ with only a mild increase as $N_s$ increases.

The performance of the solution methods is summarized in Table~\ref{tab:ex3_time_compact}.
Newton's method consistently results in the largest cumulative number of outer iterations
$K_{\mathrm{out}}$, significantly exceeding all nonlinearly preconditioned variants.
The one-level RASPEN method exhibits a rapid growth in cumulative GMRES iterations
$K_{\mathrm{gm}}$ as the number of subdomains $N_s$ increases, more than doubling when $N_s$
is doubled in one spatial direction.
In contrast, the two-level RASPEN method yields the lowest cumulative GMRES cost for the larger
partitionings, while requiring only a moderate number of additional coarse linear solves $K_c$.
Although two-level RASPEN does not achieve the same reduction in outer iterations as in the
previous two examples, it remains the most efficient method overall for large $N_s$,
particularly in terms of total fine-scale linear solver effort.

\paragraph{Time evolution and effect of flooding}
For all methods except two-level RASPEN, \Cref{fig:linsolvesperit} shows a gradual increase of $k_{\mathrm{gm}}$ over time,
which correlates with the growth of the flooded region and the resulting increase in problem difficulty.
The relative stability of two-level RASPEN suggests that applying the coarse correction nonlinearly
provides a more effective mitigation of this increasing difficulty than purely linear coarse corrections.

\paragraph{Anderson dependence on $N_s$}

\begin{figure} 
		\centering
			\includegraphics[height=6cm]{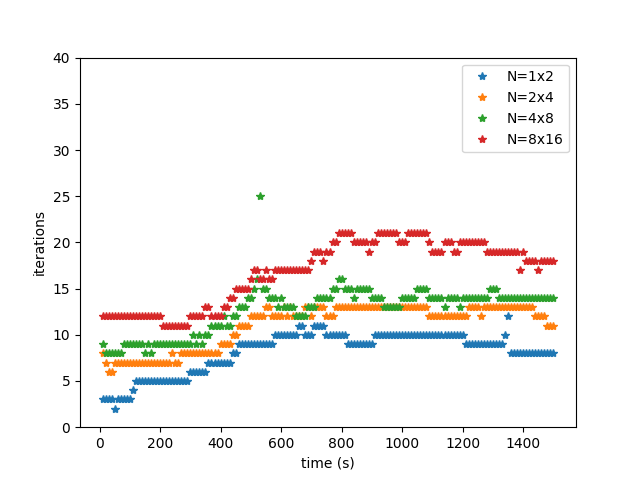}
		\caption{Test case \ref{sec:LargeUrban_num}, Anderson outer iterations per time step. Colors correspond to various coarse partitionings.}
		\label{fig:anderson}
	\end{figure}

We report in \Cref{fig:anderson} the average outer iterations per time step for Anderson method.
While the iteration count increases with $N_s$, it remains relatively stable over time.
This behavior suggests that the coarse correction used in Anderson could be further improved for this model problem.

\begin{table}[h]
\centering
\caption{Test case \ref{sec:LargeUrban_num} (time-dependent): cumulative iteration counts over the full simulation.
We report the number of time steps $N_T$, cumulative outer iterations $K_{\mathrm{out}}=N_T k_{\mathrm{out}}$,
cumulative GMRES iterations $K_{\mathrm{gm}}=N_T k_{\mathrm{out}} k_{\mathrm{gm}}$ (when applicable),
and cumulative local linear solves per subdomain $K_{\mathrm{loc}}=N_T k_{\mathrm{out}} k_{\mathrm{loc}}$.
For two-level RASPEN, the value in parentheses indicates $K_c=N_T k_{\mathrm{out}} k_c$.}
\label{tab:ex3_time_compact}
\renewcommand{\arraystretch}{1.15}
\setlength{\tabcolsep}{3pt}
\begin{tabular}{lcc|cc|cc|cc|cc}
\hline
& \multicolumn{2}{c}{$N_s=1\times2$}
& \multicolumn{2}{c}{$N_s=2\times4$}
& \multicolumn{2}{c}{$N_s=4\times8$}
& \multicolumn{2}{c}{$N_s=8\times16$} \\
Method
& $K_{\mathrm{out}}$ & $K_{\mathrm{gm}}$
& $K_{\mathrm{out}}$ & $K_{\mathrm{gm}}$
& $K_{\mathrm{out}}$ & $K_{\mathrm{gm}}$
& $K_{\mathrm{out}}$ & $K_{\mathrm{gm}}$ \\
\hline
Newton
& 1791 & 21608
& 1776 & 26753
& 1807 & 33432
& 1797 & 44146 \\
Anderson
& 1225 & --
& 1657 & --
& 1902 & --
& 2524 & -- \\
Two-step
& 597 & 7794
& 592 & 9483
& 638 & 13164
& 687 & 20029 \\
1-level RASPEN
& 477 & 3827
& 571 & 9040
& 624 & 33117
& 686 & 100453 \\
2-level RASPEN
& 515 & 3983\,(+1114)
& 595 & 6328\,(+1755)
& 650 & 8065\,(+2031)
& 730 & 12117\,(+2907) \\
\hline
\end{tabular}
\end{table}

\section{Summary and Future Work}\label{sec:conc}

This work has introduced and studied several nonlinear domain decomposition strategies
for solving nonlinear partial differential equations on complex perforated domains.
The focus has been on the Diffusive Wave model posed on urban geometries,
where perforations represent buildings and walls.
We have demonstrated that the coarse space introduced in \cite{boutilier},
originally designed for linear problems on perforated domains,
can be effectively embedded into two-level nonlinear solvers
and yields robust performance for this challenging nonlinear model.

\paragraph{GMRES iterations and scalability}
From the numerical experiments, we observe that the one-level RASPEN method
is not robust with respect to the number of subdomains $N_s$,
as evidenced by a rapid growth of the average GMRES iteration count $k_{\mathrm{gm}}$
and the cumulative GMRES cost as $N_s$ increases.
For the linearized systems arising in Newton and two-step methods,
the Trefftz coarse space performs well as a component of the two-level linear RAS preconditioner,
leading to moderate values of $k_{\mathrm{gm}}$ across all test cases.
Nevertheless, a mild increase of $k_{\mathrm{gm}}$ is observed as $N_s$ grows,
suggesting that updating or enriching the Trefftz basis as the linearization evolves
could further improve robustness.

The two-level RASPEN method consistently yields the smallest cumulative number of GMRES iterations,
particularly for large $N_s$.
Although Newton’s method exhibits comparable values of $k_{\mathrm{gm}}$ per outer iteration,
its significantly larger number of outer iterations $k_{\mathrm{out}}$
results in a substantially higher total GMRES cost.
As a consequence, Newton’s method is overall the most expensive strategy
in terms of fine-scale linear solves.

\paragraph{Outer iterations and nonlinear convergence}
In terms of outer iterations, Newton’s method is markedly more sensitive to the initial guess
than all other methods considered, with variations of up to $50\%$ in $k_{\mathrm{out}}$
across the tested configurations.
Moreover, Newton generally requires significantly more outer iterations
than the two-step and RASPEN variants.
The introduction of a nonlinear RAS correction in the two-step method
substantially accelerates convergence compared to Newton’s method,
allowing the solver to enter the fast-convergence regime much earlier.

For the stationary test cases, the two-level RASPEN method provides
the strongest reduction in $k_{\mathrm{out}}$ among all methods considered.
This effect is less pronounced for the time-dependent problem,
likely due to the influence of the time-stepping strategy
and the chosen time increment.
Nonetheless, even in the time-dependent setting,
two-level RASPEN remains the most efficient method
in terms of overall linear solver effort.

Anderson acceleration of the coarse two-step method
offers an interesting alternative to Newton’s method,
particularly for stationary problems,
where it reduces the number of outer iterations.
Each Anderson iteration involves only a small least-squares problem
and is therefore cheaper than a Newton iteration.
However, in the time-dependent example,
Anderson acceleration exhibits a stronger dependence on $N_s$
and does not provide the same level of robustness or acceleration
as the two-level RASPEN method.

\paragraph{Implementation aspects}
From an implementation standpoint,
Newton’s method and the two-step method both require the solution
of a sparse linear system at each outer iteration,
which may be solved by a direct method when problem sizes permit.
The two-step method is particularly straightforward to implement,
differing from Newton’s method only by the addition
of a nonlinear RAS correction at each iteration.
The two-level RASPEN method is more involved,
as it combines nonlinear RAS updates,
a coarse nonlinear problem,
and a globally coupled linear solve that is typically addressed by a Krylov method.
Once implemented, however,
two-level RASPEN consistently provides the most robust and efficient solver.
Anderson acceleration benefits from the availability of mature
and easy-to-use implementations in many scientific computing environments.

\paragraph{Limitations and outlook}
We note that the dimension of the coarse space increases due to edge splitting
induced by perforations, which is a known limitation of geometrically robust coarse spaces.
Alternative approaches, such as spectral or energy-minimizing coarse spaces,
may reduce the coarse dimension at the cost of increased setup complexity.
Investigating such alternatives for nonlinear problems on perforated domains
constitutes an interesting direction for future research.

Finally, we have shown that local adaptive time-step reduction
is an effective strategy for avoiding global time-step reductions
in time-dependent simulations.
Future work will focus on fully parallel implementations
and detailed wall-clock performance measurements.

\section*{Acknowledgments}

This work has been supported by ANR Project Top-up (ANR-20-CE46-0005). The high-resolution structural data has been provided by M\'{e}tropole Nice C\^{o}te d'Azur.  We warmly thank Florient Largeron, chief of MNCA's
SIG 3D project,  for his help in preparation of the data and for the multiple fruitful
discussions.

\bibliography{biblioabbrevFV}
\bibliographystyle{siam}
		
\end{document}